\documentclass[11pt]{amsart}
\usepackage[letterpaper, total={6.5in,9in}, includefoot, centering]{geometry}
\usepackage{amsfonts}
\usepackage{amsmath,amsthm,amssymb,amscd}
\usepackage[foot]{amsaddr}

\usepackage[mathscr]{eucal}
\usepackage{latexsym}
\usepackage[new]{old-arrows}
\usepackage{graphics}
\usepackage{verbatim}
\usepackage{upgreek}
\usepackage{overpic}
\usepackage[all,cmtip]{xy} 
\usepackage{enumitem}
\usepackage{tikz}
\usepackage{tikz-cd}
\usepackage[pagebackref]{hyperref}
\renewcommand*\backref[1]{\ifx#1\relax \else (Page #1) \fi}
\usepackage[font=scriptsize]{caption}
\usepackage{upgreek}
\usepackage{algorithm}
\usepackage{algpseudocode}
\makeatletter
\newsavebox{\@brx}
\newcommand{\llangle}[1][]{\savebox{\@brx}{\(\m@th{#1\langle}\)}%
  \mathopen{\copy\@brx\kern-0.5\wd\@brx\usebox{\@brx}}}
\newcommand{\rrangle}[1][]{\savebox{\@brx}{\(\m@th{#1\rangle}\)}%
  \mathclose{\copy\@brx\kern-0.5\wd\@brx\usebox{\@brx}}}
\makeatother
\graphicspath{ {./images/} }

\usepackage[authoryear,sort&compress]{natbib}
\usepackage{adjustbox}
\bibpunct{[}{]}{;}{n}{,}{,}
\usepackage{cancel}
\usepackage{cleveref}

\makeatletter \@addtoreset{equation}{section}

\makeatother
\newtheorem{theorem}{Theorem}[section]

\newtheorem{corollary}{Corollary}[section]

\newtheorem{remark}{Remark}[section]
\newtheorem{example}{Example}[section]
\newtheorem{prop}{Proposition}[section]

\begin{document}
\title{A Type II Hamiltonian Variational Principle and Adjoint Systems for Lie Groups}
\author[1]{Brian K. Tran$^1$ and Melvin Leok$^2$}
\address{$^1$Los Alamos National Laboratory, Theoretical Division, Los Alamos, NM 87545, USA. \newline \indent $^2$Department of Mathematics, University of California, San Diego, 9500 Gilman Drive, La Jolla, CA 92093-0112, USA.}
\email{btran@lanl.gov, mleok@ucsd.edu}
\allowdisplaybreaks

\begin{abstract}
We present a novel Type II variational principle on the cotangent bundle of a Lie group which enforces Type II boundary conditions, i.e., fixed initial position and final momentum. In general, such Type II variational principles are only globally defined on vector spaces or locally defined on general manifolds; however, by left translation, we are able to define this variational principle globally on cotangent bundles of Lie groups. Type II boundary conditions are particularly important for adjoint sensitivity analysis, which is our motivating application. As such, we additionally discuss adjoint systems on Lie groups, their properties, and how they can be used to solve optimization problems subject to dynamics on Lie groups. 
\end{abstract}

\maketitle


\section{Introduction}
We present a Type II Hamiltonian variational principle on the cotangent bundle of a Lie group which enforces the dynamics of a Hamiltonian system subject to initial position and fixed final momentum boundary conditions, using the trivialization of the cotangent bundle by left trivialization. For the particular choice of Hamiltonian given by the adjoint Hamiltonian associated with an ODE on a Lie group, this Type II variational principle will yield the trivialized adjoint system for dynamics on a Lie group with Type II boundary conditions. We discuss properties this system and how it naturally arises in optimization problems on Lie groups. As an example, we will show how discretizing the Type II variational principle results in Lie group variational integrators (see \cite{LeLeMc2007a, LeLeMc2007b, LeLeMc2005a, BoSu1999, MaPeSh1999, MaRo2010, KoMa2011}) and furthermore, discuss such integrators in the context of adjoint systems.

\subsection{Adjoint Systems}

The solution of many nonlinear problems involves successive linearization, and as such variational equations and their adjoints play a critical role in a variety of applications. Adjoint equations are of particular interest when the parameter space has significantly higher dimension than that of the output or objective. In particular, the simulation of adjoint equations arise in sensitivity analysis~\cite{Ca1981, CaLiPeSe2003}, adaptive mesh refinement~\cite{LiPe2003}, uncertainty quantification~\cite{WaDuAlIa2012}, automatic differentiation~\cite{Gr2003}, superconvergent functional recovery~\cite{PiGi2000}, optimal control~\cite{Ro2005}, optimal design~\cite{GiPi2000}, optimal estimation~\cite{NgGeBe2016}, and deep learning viewed as an optimal control problem~\cite{DeCeEhOwSc2019}.

The study of geometric aspects of adjoint systems arose from the observation that the combination of any system of differential equations and its adjoint equations are described by a formal Lagrangian~\cite{Ib2006, Ib2007}. This naturally leads to the question of when the formation of adjoints and discretization commutes~\cite{SoTz1997}, and prior work on this include the Ross--Fahroo lemma~\cite{RoFa2001}, and the observation by \citet{Sa2016} that the adjoints and discretization commute if and only if the discretization is symplectic.

In \cite{TrLe2022adj}, we investigated the symplectic and presymplectic geometry of adjoint systems associated with ODEs and DAEs on finite-dimensional manifolds, respectively. In \cite{TrSoLe2024}, we extended this investigation to adjoint systems associated with evolutionary PDEs on infinite-dimensional spaces. Here, we will consider adjoint systems associated with dynamics on Lie groups specifically, as they arise broadly in control and optimization. Since Lie groups are parallelizable, we will investigate adjoint systems on the trivialization of the cotangent bundle $T^*G$ of a Lie group $G$. We now introduce such systems.

\begin{example}[Adjoint Systems on Lie Groups]\label{Example2adjoint}
Consider an ODE on a Lie group $G$ given by $\dot{g} = F(g)$, specified by a vector field $F$ on $G$. We associate with $F$ the adjoint Hamiltonian $H: T^*G \rightarrow \mathbb{R}$, given by
$$ H(g,p) = \langle p, F(g)\rangle, $$
where $(g,p)$ are canonical coordinates on $T^*G$. We refer to the Hamiltonian system $i_{X_H}\Omega = dH$, relative to the canonical symplectic form $\Omega$ on $T^*G$, as the adjoint system associated with the ODE $\dot{g} = F(g)$. In coordinates, the adjoint system has the form
\begin{align*}
    \dot{g} &= F(g), \\
    \dot{p} &= -DF(g)^*p,
\end{align*}
where $DF(g)$ denotes the linearization of $F$ at $g$ and $Df(g)^*$ its adjoint with respect to the pairing $\langle\cdot,\cdot\rangle: T^*_gG \times T_gG \rightarrow \mathbb{R}$. For a time interval $[0,T]$, we will be interested in Type II boundary conditions for the above system, i.e., $g(0)=g_0$ and $p(T)=p_1$.
The motivation for considering Type II boundary conditions, arises from the fact that, viewing the ODE on $G$ as flowing forward in time, the momentum or adjoint variable flows backwards in time and can be used to extract derivative information of cost functions in the context of optimization problems subject to the dynamics on the Lie group. We will describe this in more detail in Section \ref{AdjointEquationsSection}. 
\end{example}

\section{A Type II Variational Principle on Cotangent Bundles of Lie Groups}\label{sec:type-ii-vp}
A Type II variational principle for Hamiltonian systems on the vector space $T^*Q \cong Q \times Q^*$ was introduced in \cite{LeZh2009}; however, this variational principle relies on a quantity of the form $\langle p(T), q(T)\rangle$ which does not make intrinsic sense on the cotangent bundle of a manifold $T^*M$. Furthermore, such a variational principle requires specifying a terminal condition on the covector $p(T) = p_1$ without specifying its basepoint $q(T)$, which similarly does not make intrinsic sense on $T^*M$. We present an alternative Type II variational principle which remedies these issues, in the case of the cotangent bundle of a Lie group.

To develop a variational principle for Hamiltonian systems on $T^*G$, we first consider the boundary conditions that we wish to place on the system. Note that fixed endpoint conditions on the basespace $g(0) = g_0, g(T) = g_1$ are generally incompatible with systems of the form given in Example \ref{Example2adjoint}, since adjoint systems on $T^*G$ cover first-order ODEs on $G$ and thus, one cannot freely specify both $g(0)$ and $g(T)$. As such, we instead consider Type II boundary conditions of the form $g(0) = g_0, p(T) = p_1$. As mentioned above, the issue with these boundary conditions is that one cannot intrinsically specify a covector $p(T) = p_1$ without specifying the basepoint $q(T) = q_1$. This is not an issue for adjoint systems in particular, since the time-$T$ flow of the underlying ODE on $G$ determines the basepoint where $p_1$ is specified. However, since we would like this variational principle to apply to general Hamiltonian systems on $T^*G$, we do not restrict ourselves to adjoint systems in this section. Fortunately, we can make sense of Type II boundary conditions on Lie groups, since $T^*G$ is trivializable by left-translation. 

Let $\mathfrak{g} = T_eG$ denote the Lie algebra of $G$ and $\mathfrak{g}^* = T^*_eG$ be its dual. We will denote the duality pairing between $v \in T_gG$ and $p \in T_g^*G$ by $\langle p,v\rangle$, where the base point is understood in context. Let $L_g: G \rightarrow G$ denote left-translation by $g$,  $L_g(x)=gx$. Left-translation induces maps on the tangent bundle and cotangent bundle of $G$ by pushforward and pullback, respectively, which we denote by
\begin{align*}
T_xL_g&: T_xG \rightarrow T_{gx}G, \\
T_x^*L_g&: T_{gx}^*G \rightarrow T_x^*G.
\end{align*}
For $v_g \in T_gG, p_g \in T_g^*G$, we will adopt the following compact notation for their left-translations to their respective fibers over the identity,
\begin{align*}
g^{-1} v_g &\equiv T_gL_{g^{-1}}(v_g) \in \mathfrak{g}, \\
g^* p_g &\equiv T^*_eL_g(p_g) \in \mathfrak{g}^*.
\end{align*}
This notation is suggestive, since in the case that $G$ is a matrix Lie group, the left-translation of a tangent vector to the fiber over the identity acts by matrix multiplication by the inverse of $g$ and the left-translation of a covector to the fiber over the identity acts by matrix multiplication by the adjoint of $g$.

A useful fact is that the pairing $\langle p_g, v_g\rangle$ is preserved under left-translation,
\begin{align*}
\langle g^* \cdot p_g, g^{-1}\cdot v_g\rangle &= \langle T_e^*L_g (p_g), T_gL_{g^{-1}}(v_g)\rangle = \langle p_g, T_eL_g \circ T_gL_{g^{-1}} (v_g)\rangle \\ 
&= \langle p_g, T_g(L_g \circ L_{g^{-1}})(v_g)\rangle = \langle p_g, v_g\rangle. 
\end{align*}
By left-translation on the cotangent bundle, we get the trivialization $T^*G \cong G \times \mathfrak{g}^*$. With this trivialization, we can make sense of Type II boundary conditions $g(0) = g_0 \in G, \mu(T) = \mu_1 \in \mathfrak{g}^*$, with coordinates $(g,\mu)$ on $G \times \mathfrak{g}^*$.

What remains is to construct a variational principle. Recall that the action for a Hamiltonian system on $T^*G$ is given by
$$ S[g,p] = \int_0^T \Big( \langle p, \dot{g}\rangle - H(g,p) \Big) dt, $$
where $H: T^*G \rightarrow \mathbb{R}$. By left-translation, with $\mu = g^*\cdot p$, we define the left-trivialized Hamiltonian $h: G \times \mathfrak{g}^* \rightarrow \mathbb{R}$ as 
$$ h(g,\mu) \equiv H(g,g^{*-1}\cdot \mu) = H(g,p). $$
The action can then be expressed as
\begin{align*}
S[g,p] &= \int_0^T \Big( \langle p, \dot{g}\rangle - H(g,p) \Big) dt = \int_0^T \Big( \langle g^*\cdot p, g^{-1}\cdot \dot{g}\rangle - H(g,p) \Big) dt \\
&= \int_0^T \Big( \langle \mu, g^{-1}\cdot \dot{g}\rangle - h(g,\mu) \Big) dt =: s[g,\mu].
\end{align*}
We refer to $s[g,\mu]$ as the left-trivialized action. 

Now, we prescribe boundary conditions $g(0) = g_0 \in G$, $\mu(T) = \mu_1 \in \mathfrak{g}^*$. Given a curve $(g(t),p(t))$ on $T^*G$, by left-translation, the terminal momenta condition $\mu(T) = \mu_1$ on $\mathfrak{g}^*$ corresponds to $p(T) = g(T)^{*-1}\cdot \mu_1 \in T^*_{g(T)}G$. To state a variational principle, we observe that by left-translation, we can prescribe a boundary condition on $\mu(T)$ which, under the inverse of the trivialization, corresponds to a boundary condition on $p(T)$. Now, consider variations $(\delta g, \delta p)$ of a curve $(g,p)$ satisfying $\delta g(0) = 0$ and $\delta p(T) = 0$, which correspond to Type II boundary conditions. Since $\delta g(T)$ is arbitrary (as opposed to traditional variational principles which fix both position endpoints), virtual work can always be performed on this system through such variations, i.e., the variation of the action under such variations would be $\langle p(T), \delta g(T)\rangle$, or equivalently, $\langle \mu(T), \eta(T)\rangle$ where we defined the left-trivialization of the variation $\eta = g^{-1}\cdot \delta g$. Thus, we cannot demand that the action $S$ (equivalently, $s$) is stationary since one can always introduce virtual work as described above; however, we can demand that it is stationary modulo the virtual work that is introduced into the system by varying the terminal point $g(T)$. Thus, we impose the variational principle
$$ \delta S[g,p] = \langle p(T),\delta g(T)\rangle, $$
or equivalently, by left translation
$$ \delta s[g,\mu] = \langle \mu(T), \eta(T)\rangle, $$
subject to variations fixing $g(0)$ and $p(T)$ (equivalently, $\mu(T)$). We refer to this variational principle as the Type II d'Alembert variational principle, due to its similarity to the d'Alembert variational principle which utilizes virtual work to derive forced Lagrangian or Hamiltonian systems~\cite{MaWe2001}.

\begin{theorem}[Type II d'Alembert Variational Principle]\label{TypeIIdAVP Theorem}
The following are equivalent:
\begin{enumerate}[label=\normalfont(\roman*)]
\item The Type II d'Alembert variational principle 
$$\delta S[g,p] = \langle p(T), \delta g(T)\rangle$$
on $T^*G$ is satisfied, where the variations $(\delta g, \delta p)$ of the curve $(g,p)$ satisfy $\delta g(0) = 0, \delta p(T) = 0$, which correspond to boundary conditions $g(0)=g_0, p(T) = g(T)^{*-1}\cdot \mu_1$.
\item Hamilton's equations hold in canonical coordinates on $T^*G$, with the above Type II boundary conditions,
\begin{subequations}
\begin{align}
\dot{g} &= D_p H(g,p), \label{TypeIIHamilton a}  \\
\dot{p} &= - D_g H(g,p), \label{TypeIIHamilton b} \\
g(0) &= g_0, \label{TypeIIHamilton c}\\
p(T) &= g(T)^{*-1}\cdot \mu_1. \label{TypeIIHamilton d}
\end{align}
\end{subequations}
\item The Type II d'Alembert variational principle
$$ \delta s[g,\mu] = \langle \mu(T), \eta (T)\rangle $$
on $G \times \mathfrak{g}^*$ is satisfied, where the variation $\delta g$ is left-trivialized as $\eta = g^{-1}\cdot \delta g$ and the variation $\delta p$ is left-trivialized as $\delta \mu = g^* \cdot \delta p$, with $\delta \eta (0)=0, \delta \mu(T) = 0$, which correspond to boundary conditions $g(0)=g_0, \mu(T) = \mu_1$. 
\item The Lie--Poisson equations hold on $G \times \mathfrak{g}^*$, with the above Type II boundary conditions,
\begin{subequations}
\begin{align}
\dot{g} &= g \cdot D_\mu h(g,\mu), \label{TypeIILPa} \\
\dot{\mu} &= -g^* \cdot D_g h(g,\mu) + \textup{ad}^*_{D_\mu h(g,\mu)}\mu, \label{TypeIILPb} \\
g(0) &= g_0, \label{TypeIILPc} \\
\mu(T) &= \mu_1. \label{TypeIILPd}
\end{align}
\end{subequations}
\end{enumerate}
\begin{remark}
Above, we denote by $D_gH, D_pH, D_gh, D_\mu h$ the derivatives of $H$, which are related to the exterior derivative $dH$ of $H$ by
\begin{align*}
dH(g,p)\cdot (\delta g,\delta p) = \langle D_gH(g,p), \delta g\rangle + \langle \delta p, D_pH(g,p)\rangle, \\
dh(g,\mu)\cdot (\delta g,\delta \mu) = \langle D_gh(g,\mu), \delta g\rangle + \langle \delta \mu, D_\mu h(g,\mu)\rangle.
\end{align*}
\end{remark}

\begin{proof}
To see that (i) and (ii) are equivalent, compute the variation of $S$,
\begin{align*}
\delta S[g,p] &= \int_0^T \Big( \langle \delta p, \dot{g}\rangle + \left\langle p, \frac{d}{dt} \delta g\right\rangle - \langle D_gH, \delta g\rangle - \langle \delta p, D_pH\rangle \Big) dt \\
&= \int_0^T \Big( \langle \delta p, \dot{g} - D_pH\rangle + \langle - \dot{p} - D_gH, \delta g\rangle\Big) dt + \langle p,\delta g\rangle\Big|_0^T.
\end{align*}
If (ii) holds, the integrand above vanishes by the equations of motion; furthermore, $\delta g(0) = 0$. Thus, the above expression reduces to
$$ \delta S[g,p] = \langle p(T), \delta g(T)\rangle, $$
i.e., (i) holds. Conversely, if (i) holds, we have
$$ 0 = \delta S[g,p] - \langle p(T), \delta g(T)\rangle = \int_0^T \Big( \langle \delta p, \dot{g} - D_pH\rangle + \langle - \dot{p} - D_gH, \delta g\rangle\Big) dt. $$
Then, by the fundamental lemma of the calculus of variations, (ii) holds. 

To see that (iii) is equivalent to (iv), compute the variation of $s$. For simplicity, we denote the left-translation of $\dot{g}$ by $\xi = g^{-1}\cdot \dot{g}$ and similarly $\eta = g^{-1}\cdot \delta g$.
\begin{align*}
\delta s[g,\mu] &= \int_0^T \Big( \langle \delta \mu, g^{-1}\cdot \dot{g}\rangle + \left\langle \mu, g^{-1} \frac{d}{dt} \delta g - g^{-1} \cdot \delta g g^{-1}\cdot \dot{g} \right\rangle - \langle \delta \mu, D_\mu h\rangle  - \langle D_gh, \delta g\rangle \Big) dt \\
&= \int_0^T \Big( \langle \delta \mu, g^{-1}\cdot \dot{g} - D_\mu h\rangle + \langle \mu, \dot{\eta} + \text{ad}_\xi \eta \rangle - \langle g^*\cdot D_gh, \eta\rangle \Big)dt \\
&= \int_0^T \Big( \langle \delta \mu, g^{-1}\cdot \dot{g} - D_\mu h\rangle + \langle -\dot{\mu} + \text{ad}_{\xi}^*\mu - g^*\cdot D_g h, \eta\rangle \Big)dt + \langle \mu,\eta\rangle\Big|_0^T.
\end{align*}
If (iv) holds, the integrand above vanishes by the equations of motion, noting that $\xi = g^{-1}\cdot \dot{g} = D_\mu h$; furthermore, $\eta(0)=0$. Thus, the above expression reduces to
$$ \delta s[g,\mu] = \langle \mu(T),\eta(T)\rangle, $$
i.e., (iii) holds. Conversely, if (iii) holds, we have
$$ 0 = \delta s[g,\mu] - \langle \mu(T),\eta(T)\rangle = \int_0^T \Big( \langle \delta \mu, g^{-1}\cdot \dot{g} - D_\mu h\rangle + \langle -\dot{\mu} + \text{ad}_{\xi}^*\mu - g^*\cdot D_g h, \eta\rangle \Big)dt. $$
Then, by the fundamental lemma of the calculus of variations, (iv) holds.

Finally, (i) and (iii) are equivalent by left-translation, since $S[g,p]=s[g,\mu]$ and $\langle p(T),\delta g(T)\rangle = \langle \mu(T),\eta(T)\rangle$.
\end{proof}
\end{theorem}

\begin{remark}
Note that one can also modify the above variational principle to include external forces by adding the virtual work done by the external force. Given a left-trivialized external force $f: [0,T]\rightarrow \mathfrak{g}^*$, one can modify the above variational principle to 
$$ \delta s[g,\mu] = \langle \mu(T), \eta(T)\rangle + \int_0^T \langle f, \eta\rangle dt, $$
or equivalently, 
$$ \delta S[g,p] = \langle p(T), \delta g(T)\rangle + \int_0^T \langle g^{*-1}\cdot f, \delta g\rangle dt. $$
This modifies the momenta equations \eqref{TypeIILPb} on $G \times \mathfrak{g}$ to include the external force,
$$ \dot{\mu} = -g^* \cdot D_g h(g,\mu) + \textup{ad}^*_{D_\mu h(g,\mu)}\mu + f, $$
or equivalently, modifies the momenta equation \eqref{TypeIIHamilton b} on $T^*G$ to be
$$ \dot{p} = - D_g H(g,p) + g^{*-1} \cdot f.$$
\end{remark}

A particularly important class of Hamiltonians are the left-invariant Hamiltonians, which are functions $H: T^*G \rightarrow \mathbb{R}$ that are invariant under the cotangent lift of left-multiplication by any $x \in G$, i.e.,
$$ H \circ T^*L_x = H \text{ for all } x \in G. $$
In terms of our notation, that is
$$ H(xg,x^{*-1}p)=H(g,p) \text{ for all } x \in G,\ (g,p) \in T^*G. $$
For such left-invariant Hamiltonians, the dynamics on $T^*G$ reduce to dynamics on $\mathfrak{g}^*$ \cite{MaRa1999}.

Given a left-invariant Hamiltonian $H$, we define the reduced Hamiltonian $\tilde{H}: \mathfrak{g}^* \rightarrow \mathbb{R}$ by
$$ \tilde{H}(\mu) = H(e,\mu). $$
Then, equation \eqref{TypeIILPb} for $\dot{\mu}$ reduces to 
\begin{equation}\label{Continuous Reduced Momentum Equation}
\dot{\mu} = \text{ad}^*_{D_\mu \tilde{H}(\mu)} \mu,
\end{equation}
where we used that 
$$\tilde{H}(\mu) = H(e,\mu) = H(g, g^{*-1}\mu) = h(g,\mu),$$
and hence, $D_gh(g,\mu) = 0$. Thus, as can be seen from equation \eqref{Continuous Reduced Momentum Equation}, the momentum equation decouples from the dynamics on $G$: \eqref{Continuous Reduced Momentum Equation} with $\mu(T) = \mu_1$ can be solved independently and subsequently, \eqref{TypeIILPa} can be used to reconstruct the dynamics on $G$. Hence, for a left-invariant system, the full dynamics on $T^*G$ are completely encoded by the reduced dynamics on $\mathfrak{g}^*$.

\section{Adjoint Systems on Lie Groups}\label{AdjointEquationsSection}
In this section, we develop the theory of adjoint systems and sensitivity analysis on Lie groups. We thus focus on the case where the Hamiltonian system on $T^*G$ is an adjoint system, as introduced in Example \ref{Example2adjoint}.

Let $F$ be a vector field on $G$ and consider the differential equation $\dot{g} = F(g)$ on $G$. We define the adjoint Hamiltonian associated with $F$ as
\begin{align*}
H: T^*G &\rightarrow \mathbb{R}, \\
  (g,p) & \mapsto H(g,p) \equiv \langle p, F(g)\rangle. 
\end{align*}
In canonical coordinates $(g,p)$ on $T^*G$, the adjoint system \eqref{TypeIIHamilton a}-\eqref{TypeIIHamilton b} has the form
\begin{align*}
\dot{g} &= F(g), \\
\dot{p} &= - DF(g)^* p.
\end{align*}

We begin by computing the Lie--Poisson equations \eqref{TypeIILPa}-\eqref{TypeIILPb} for this particular class of adjoint Hamiltonians. We denote by $f$ the left-trivialization of $F$,
\begin{align*}
f: G &\rightarrow \mathfrak{g}, \\
g &\mapsto f(g)\equiv g^{-1} \cdot F(g).
\end{align*}
Then, the left-trivialized Hamiltonian $h: G \times \mathfrak{g}^* \rightarrow \mathbb{R}$ is given by $h(g,\mu) = \langle \mu, f(g)\rangle.$ Computing the functional derivatives of $h$ yields
\begin{align*}
D_\mu h(g,\mu) &= f(g), \\
D_g h(g,\mu) &= Df(g)^*\mu.
\end{align*}
Thus, the Lie--Poisson system \eqref{TypeIILPa}-\eqref{TypeIILPb} for the adjoint Hamiltonian has the form
\begin{subequations}
\begin{align}
\dot{g} &= F(g), \label{LP Adjoint a}\\
\dot{\mu} &= - g^*\cdot [Df(g)]^* \mu + \text{ad}^*_{f(g)}\mu. \label{LP Adjoint b}
\end{align}
\end{subequations}

We now address the question of existence and uniqueness for solutions of the Type II system \eqref{TypeIILPa}-\eqref{TypeIILPd}. For general Hamiltonians on $T^*G$, this is a complicated question which is dependent on the particular Hamiltonian. In particular, since the system has Type II boundary conditions $g(0)=g_0, \mu(T)=\mu_1$, it can be interpreted as a two-point boundary value problem whose solution theory is in general problem dependent (see, e.g., \cite{BaSh1969, LaOp1967, ElHe2016}), as opposed to initial-value problems with $g(0) = g_0, \mu(0)=\mu_0$ where the solution would simply be given by the flow of the vector field.

However, for adjoint systems in particular, we can provide a global solution theory which utilizes the fact that the adjoint system covers an ODE on $G$. Assuming the ODE on $G$ behaves nicely, we will have unique solutions for the adjoint system on $T^*G$. We make this more precise in the following proposition.

\begin{prop}[Global Existence and Uniqueness of Solutions to Adjoint Systems on $T^*G$]\label{Solution to Adjoint System Proposition}

Let $T>0, g_0 \in G, \mu_1 \in \mathfrak{g}^*$. Let $F$ be a complete vector field on $G$, i.e., it generates a global flow $\Phi_F: \mathbb{R} \times G \rightarrow G$.

Then, there exists a unique curve $(g,\mu): [0,T] \rightarrow G \times \mathfrak{g}^*$ satisfying the Lie--Poisson system with Type II boundary conditions \eqref{TypeIILPa}-\eqref{TypeIILPd}, where $h$ is the left-trivialized adjoint Hamiltonian associated with $F$. 

Furthermore, there exists a unique curve $(g,p): [0,T] \rightarrow T^*G$ satisfying Hamilton's equations with Type II boundary conditions \eqref{TypeIIHamilton a}-\eqref{TypeIIHamilton d}, where $H$ is the adjoint Hamiltonian associated with $F$.
\begin{proof}
By the fundamental theorem on flows \cite{Le2012}, there exists a unique curve $g: \mathbb{R} \rightarrow G$ satisfying $\dot{g} = F(g)$ and $g(0) = g_0$, given by the flow of $F$ on $g_0$, $g(t) = \Phi_F(t,g_0)$. In particular, $g$ is a smooth function of $t$, since $F$ is smooth. Recall that we assume all maps and manifolds are smooth, unless otherwise stated.

Now, with this curve $g(t)$ fixed, we substitute this into the differential equation for $\mu$ \eqref{TypeIILPb}, to obtain
$$ \dot{\mu} = -g(t)^*\cdot [Df(g(t))]^* \mu + \text{ad}^*_{f(g(t))}\mu. $$
In particular, this equation has the form of a time-dependent linear differential equation on $\mathfrak{g}^*$,
$$ \dot{\mu} = L(t)\mu, $$
where we define the time-dependent linear operator $L: \mathbb{R} \rightarrow \text{End}(\mathfrak{g}^*)$ by
\begin{equation}\label{linear_operator}
	L(t) = -g(t)^*\cdot [Df(g(t))]^* + \text{ad}^*_{f(g(t))}.
\end{equation}
Since $g$ is a smooth function of $t$, $L$ is a smooth, and in particular continuous, function of $t$. Hence, by the standard solution theory for linear differential equations, there exists a unique curve $\mu: [0,T] \rightarrow \mathfrak{g}^*$ satisfying $\dot{\mu}= L(t)\mu$ and $\mu(T) = \mu_1$.

For the second statement of the proposition, note that solution curves $(g,p): [0,T] \rightarrow T^*G$ of \eqref{TypeIIHamilton a}-\eqref{TypeIIHamilton d} are in one-to-one correspondence with solution curves $(g,\mu): [0,T] \rightarrow G \times \mathfrak{g}^*$ of \eqref{TypeIILPa}-\eqref{TypeIILPd} via left-translation. 
\end{proof}
\end{prop}

By the above proposition, we know that there exists a unique solution to the adjoint system on $T^*G$ with Type II boundary conditions, under the assumption that $F$ is complete. For Lie groups, there are two particularly important cases where this assumption is satisfied.

\begin{corollary}
If $G$ is a compact Lie group, then the above proposition holds for any vector field $F$ on $G$.

If $F$ is a left-invariant vector field on a (not necessarily compact) Lie group $G$, then the above proposition holds.
\begin{proof}
The first statement follows from the fact that any vector field on a compact manifold is complete. The second statement follows from the fact that any left-invariant vector field on a Lie group is complete. See \cite{Le2012}.
\end{proof}
\end{corollary}

\textbf{The Variational System.}
An important property of adjoint systems is that they satisfy a quadratic conservation law, which is at the heart of the method of adjoint sensitivity analysis \cite{Sa2016}. 

To state this conservation law, we introduce the variational equation associated with an ODE $\dot{g} = F(g)$ on a Lie group $G$, which is defined to be the linearization of the ODE,
$$ \frac{d}{dt} \delta g = DF(g)\delta g. $$
We refer to the combined system
\begin{subequations}
\begin{align}
\frac{d}{dt} g &= F(g), \label{VariationalEquation a} \\
\frac{d}{dt} \delta g &= DF(g)\delta g, \label{VariationalEquation b}
\end{align}
\end{subequations}
as the variational system, which is interpreted as an ODE on $TG$.

As with the adjoint system, it will be useful to left-trivialize this system, which will give an ODE on $G \times \mathfrak{g}$. As before, let $f(g) = g^{-1}\cdot F(g)$ be the left-trivialization of $F$. Let $\eta = g^{-1}\cdot \delta g$ and let $\xi = g^{-1}\cdot \dot{g}$. As is well-known (see, for example, \cite{MaRa1999}), we have the relation
$$ \dot{\eta} = \delta \xi - [\xi,\eta]. $$
In particular, since $\xi = g^{-1} \cdot \dot{g} = f(g)$, we have $\delta \xi = Df(g)\delta g = Df(g) g \cdot \eta$, so that the above relation becomes
$$ \dot{\eta} = Df(g) g \cdot \eta - [f(g),\eta], $$
which we refer to as the left-trivialized variational equation. We refer to the combined system
\begin{subequations}
\begin{align}
\dot{g} &= F(g), \label{VariationalEquation Left-Trivial a} \\
\dot{\eta} &= Df(g) g \cdot \eta  - \text{ad}_{f(g)}\eta, \label{VariationalEquation Left-Trivial b}
\end{align}
\end{subequations}
as the left-trivialized variational system on $G \times \mathfrak{g}$. Analogous to Proposition \ref{Solution to Adjoint System Proposition} on the existence and uniqueness of solutions for adjoint systems, we have the following result. 

\begin{prop}[Global Existence and Uniqueness of Solutions to Variational Systems on $TG$]\label{Solution to Variational System Proposition}

Let $T>0, g_0 \in G, \eta_0 \in \mathfrak{g}$. Let $F$ be a complete vector field on $G$, i.e., it generates a global flow $\Phi_F: \mathbb{R} \times G \rightarrow G$.

Then, there exists a unique curve $(g,\eta): [0,T] \rightarrow G \times \mathfrak{g}$ satisfying the left-trivialized variational system \eqref{VariationalEquation Left-Trivial a}-\eqref{VariationalEquation Left-Trivial b} with initial conditions $g(0)=g_0, \eta(0) = \eta_0$.

Furthermore, there exists a unique curve $(g,\delta g): [0,T] \rightarrow TG$ satisfying the variational system \eqref{VariationalEquation a}-\eqref{VariationalEquation b} with initial conditions $g(0)=g_0, \delta g(0) = g_0 \cdot \eta_0$. 
\begin{proof}
The proof is almost identical to the proof of Proposition \ref{Solution to Adjoint System Proposition}.
\end{proof}
\end{prop}

We can now state the quadratic conservation law enjoyed by solutions of the adjoint and variational systems.

\begin{theorem}\label{Quadratic Conservation Law Theorem}
Let $(g,\mu)$ be a solution curve of the left-trivialized adjoint system and let $(g,\eta)$ be a solution curve of the left-trivialized variational system, both covering the same base curve $g$. Let $(g,p)$ and $(g,\delta g)$ be the respective solution curves for the adjoint system and variational system obtained by left-translation. Then,
\begin{subequations}
\begin{align}
\frac{d}{dt} \langle \mu(t), \eta (t) \rangle &= 0, \label{QuadraticAdjointConservationLaw a} \\
\frac{d}{dt} \langle p(t), \delta g(t)\rangle &= 0. \label{QuadraticAdjointConservationLaw b}
\end{align}
\end{subequations}
\begin{proof}
Note that it suffices to prove either \eqref{QuadraticAdjointConservationLaw a} or \eqref{QuadraticAdjointConservationLaw b}, since the duality pairing is invariant under left-translation,
$$ \langle \mu(t), \eta(t)\rangle = \langle g(t)^{*-1}\cdot p(t), g(t)\cdot \delta g(t)\rangle = \langle p(t), \delta g(t)\rangle. $$

We will prove \eqref{QuadraticAdjointConservationLaw a}. Compute
\begin{align*}
\frac{d}{dt} \langle& \mu(t), \eta (t)\rangle\\
& = \langle \dot{\mu}(t), \eta(t)\rangle + \langle \mu(t), \dot{\eta}(t)\rangle \\
&= \langle -g(t)^*\cdot [Df(g(t))]^*\mu (t) + \text{ad}^*_{f(g(t))}\mu(t), \eta(t)\rangle + \langle \mu(t), Df(g(t)) g(t)\cdot \eta(t) - \text{ad}_{f(g(t))}\eta(t)\rangle \\
&= -\langle \mu(t), Df(g(t)) g(t)\cdot \eta(t)\rangle + \langle \mu(t), Df(g(t)) g(t)\cdot \eta(t)\rangle  \\
& \qquad + \langle \mu(t), \text{ad}_{f(g(t))}\eta(t)\rangle - \langle \mu(t), \text{ad}_{f(g(t))}\eta(t)\rangle \\
&= 0. \qedhere
\end{align*}
\end{proof}
\end{theorem}

In particular, we have the following corollary of Propositions \ref{Solution to Adjoint System Proposition} and \ref{Solution to Variational System Proposition} and Theorem \ref{Quadratic Conservation Law Theorem}.

\begin{corollary}\label{corollary-adjoint}
Let $T>0, g_0 \in G, \mu_1 \in \mathfrak{g}^*, \eta_0 \in \mathfrak{g}$. Let $F$ be a complete vector field on $G$. Let $g$ be the integral curve of $F$ with initial condition $g(0)=g_0$. Then, the solution curves of the adjoint and variational systems from Propositions \ref{Solution to Adjoint System Proposition} and \ref{Solution to Variational System Proposition} satisfy the quadratic conservation law
$$ \langle \mu(0), \eta_0\rangle = \langle \mu_1, \eta(T)\rangle, $$
or, equivalently, $\langle p(0), \delta g_0\rangle = \langle p_1, \delta g(T)\rangle$ where $\delta g_0 = g_0 \eta_0$ and $p_1 = g(T)^{*-1}\mu_1$.
\end{corollary}

As we will see in Section \ref{Adjoint Sensitivity Analysis Section}, this conservation law will be the basis for adjoint sensitivity analysis on Lie groups.

\textbf{Reduction of Adjoint Systems for Left-invariant Vector Fields.} In practice, many interesting mechanical systems arise from the flow of left-invariant vector fields on Lie groups. As such, we will consider adjoint systems in the particular case where the vector field is left-invariant. First, we will show that left-invariant vector fields are in one-to-one correspondence with left-invariant adjoint Hamiltonians. Subsequently, we will state the adjoint equations in this particular case.

\begin{prop}
Let $F$ be a vector field on $G$. Then the associated adjoint Hamiltonian $H(g,p) = \langle p,F(g)\rangle$ is left-invariant if and only if $F$ is left-invariant.
\begin{proof}
Assume that $F$ is left-invariant, i.e., $F(xg) = x F(g)$ for all $x,g \in G$. Then, for any $x,g \in G, p \in T^*_gG$,
\begin{align*}
H(xg, x^{*-1}p) = \langle x^{*-1}p, F(xg)\rangle = \langle x^{*-1}p, xF(g)\rangle = \langle x^* x^{*-1}p, F(g)\rangle = \langle p,F(g)\rangle = H(g,p),
\end{align*}
i.e., $H$ is left-invariant.

Conversely, assume that $H$ is left-invariant, i.e., $H(g,p) = H(xg, x^{*-1}p)$ for all $x,g \in G, p \in T^*_gG$. Then, for any $x,g \in G, p \in T^*_gG$,
$$ \langle p,F(g)\rangle = H(g,p) = H(xg,x^{*-1}p) = \langle x^{*-1}p, F(xg)\rangle = \langle p, x^{-1} F(xg)\rangle. $$
Since $p \in T^*_gG$ is arbitrary, we have for all $x,g \in G$,
$$ F(g) = x^{-1}F(xg), $$
i.e., $xF(g) = F(xg)$, so $F$ is left-invariant.
\end{proof}
\end{prop}

Since a left-invariant vector field corresponds to a left-invariant adjoint Hamiltonian, the reduction theory discussed in Section \ref{sec:type-ii-vp} applies. Thus, the adjoint equation for the momenta $\mu$, from equation \eqref{Continuous Reduced Momentum Equation}, is given by
$$ \dot{\mu} = \text{ad}^*_{F(e)}\mu, $$
since $\tilde{H}(\mu) = H(e,\mu) = \langle \mu, F(e)\rangle$ and hence, $D_\mu \tilde{H}(\mu) = F(e)$.

\subsection{Adjoint Sensitivity Analysis on Lie Groups}\label{Adjoint Sensitivity Analysis Section}
We now consider the application of the adjoint system to the following optimization problems. We consider an initial condition optimization problem
\begin{align}\label{Adjoint Sensitivity Problem}
\min_{g_0 \in G} C(&g(T)),  \\
\text{such that } &\dot{g}(t) = F(g(t)),\ t \in (0,T), \nonumber \\
&g(0) = g_0, \nonumber
\end{align}
and an optimal control problem,
\begin{align}\label{Parameter Sensitivity Problem}
\min_{u \in U} C(&g(T)),  \\
\text{such that } &\dot{g}(t) = F(g(t),u),\ t \in (0,T), \nonumber \\
&g(0) = g_0, \nonumber
\end{align}
where in \eqref{Parameter Sensitivity Problem}, we have introduced a parameter-dependent vector field $F$.

\textbf{Initial Condition Sensitivity.} We begin with problem \eqref{Adjoint Sensitivity Problem}. We refer to $C: G \rightarrow \mathbb{R}$ as the terminal cost function. Thus, the problem \eqref{Adjoint Sensitivity Problem} is to find an initial condition $g_0 \in G$ which minimizes the cost function at the terminal-value $g(T)$, subject to the dynamics of the ODE $\dot{g} = F(g)$.

For gradient-based algorithms, one needs the derivative of the terminal cost function $C(g(T))$ with respect to the initial condition $g_0$, where we regard the terminal cost as an implicit function of $g_0$; we refer to this derivative as the \textit{initial condition sensitivity}. To compute this derivative, we will use the adjoint system and particularly, \Cref{corollary-adjoint}. We will derive the left-trivialization of the derivative of the cost function with respect to $g_0$; we will discuss in \Cref{DHVI section} how this left-trivialized derivative can be used to construct a numerical method for solving the problem \eqref{Adjoint Sensitivity Problem}.

Let $(g,p)$ be the solution to the adjoint system with $g(0)=g_0$, $p(T) = dC(g(T))$ and let $(g,\delta g)$ be the solution of the variational equation with $g(0)=g_0$, $\delta g(0) = \delta g_0$. Regarding $C(g(T))$ as an implicit function of $g_0$, its derivative in the direction $\delta g_0$ is given by
\begin{align*}
    \langle D_{g_0} C(g(T)), \delta g_0 \rangle = \left\langle dC(g(T)), \frac{\delta g(T)}{\delta g_0} \delta g_0 \right\rangle = \left\langle dC(g(T)), \delta g(T) \right\rangle,
\end{align*}
where in the first equality we used the chain rule and the fact that $(\delta g(T)/\delta g_0) \delta g_0 = \delta g(T)$ for a solution $\delta g$ of the variational equation with $\delta g(0)=\delta g_0$ which follows from the fact that the variational equation is linear in $\delta g$. Now, since $p$ solves the adjoint equation with $p(T) = dC(g(T))$, from \Cref{corollary-adjoint} the right hand side of the above equation is $\langle p(0), \delta g(0)\rangle$. That is, we have
\begin{align*}
    \langle D_{g_0} C(g(T)), \delta g_0 \rangle = \left\langle p(0), \delta g_0 \right\rangle.
\end{align*}
Since $\delta g_0$ is arbitrary, we obtain the initial condition sensitivity $p(0) = D_{g_0} C(g(T))$ and in particular, the left-trivialized derivative of the cost function with respect to $g_0$ is $\mu(0) = g_0^* D_{g_0} C(g(T))$.

\textbf{Parameter Sensitivity.} Analogously, we consider problem \eqref{Parameter Sensitivity Problem} and determining the appropriate gradient for this problem. The problem \eqref{Parameter Sensitivity Problem} is to find a parameter $u \in U$ which minimizes the terminal cost function $C(g(T))$, subject to the dynamics of the parameter-dependent ODE $\dot{g} = F(g,u)$ with $g(0)=g_0$ fixed. We assume that $F$ is continuously differentiable with respect to $u$.

For a gradient-based algorithm, one needs the derivative of the terminal cost function $C(g(T))$ with respect to the parameter $u$; we refer to this derivative as the \textit{parameter sensitivity}. We show how this can be computed using the adjoint system. Define the parameter-dependent action as
$$ S[g,p;u] \equiv \int_{0}^{T} \langle p,\dot{g} - F(g,u)\rangle dt.$$
Consider the augmented cost function, given by subtracting the parameter-dependent action from the terminal cost function,
\begin{align*}
J & \equiv C(g(T)) - S[g,p;u] = C(g(T)) - \int_{0}^{T} \langle p, \dot{g} - F(g,u) \rangle dt.
\end{align*}
By left-trivialization, this is equivalent to
\begin{equation}\label{Left-trivialized Parameter Sensitivity}
J = C(g(T)) - \int_{0}^{T} \langle \mu, \xi - f(g,u) \rangle dt,
\end{equation}
where $\xi = g^{-1}\dot{g}$ and $f$ is the left-trivialization of $F$. Observe that since the integral of \eqref{Left-trivialized Parameter Sensitivity} vanishes when $\xi = f(g,u)$, we have that the derivative of $J$ with respect to $u$ equals the derivative of $C$ with respect to $u$, subject to the variational equations, where the variation $\delta_u g$ is given by the variation of $g$ induced by varying $u$. Thus,
$$ \frac{d}{du}C(g(T)) = \frac{d}{du}J.$$
\begin{prop}
The parameter sensitivity is given by
$$ \frac{d}{du} C(g(T)) = \frac{d}{du} J = \int_{0}^{T} \langle \mu, \frac{\partial}{\partial u}f(g,u)\rangle dt, $$
where $\mu$ is chosen to satisfy the adjoint equation $ -\dot{\mu} + \text{ad}^*_{\xi}\mu - g^{*}[D_g f(g,u)]^*\mu  = 0$ and the terminal condition $\mu(T) = g(T)^* dC(g(T))$.
\begin{proof}
We compute $dJ/du$ explicitly,
\begin{align*}
\frac{d}{du}J(g(T)) &= \langle dC(g(T)), \delta_u g (T)\rangle - \int_{0}^{T} \Big[ \langle \mu, \delta_u \xi - [D_g f(g,u)] g \eta_u \rangle - \langle \mu, \frac{\partial}{\partial u} f(g,u)\rangle \Big] dt,
\end{align*}
where we introduced the left-trivialized variation $\eta_u = g^{-1}\delta_u g$ and we have decomposed the total derivative of $f$ with respect to $u$ into its implicit dependence on $u$ through $g$ as well as its explicit dependence on $u$, i.e.,
$$ \frac{d}{du} f(g,u) = [D_g f(g,u)] \delta_u g + \frac{\partial}{\partial u}f(g,u). $$
Using the relation $\dot{\eta}_u = \delta_u \xi - \text{ad}_{\xi} \eta_u$, this becomes
\begin{align*}
\frac{d}{du}J(g(T)) &= \langle dC(g(T)), \delta_u g (T)\rangle - \int_{0}^{T} \Big[ \langle \mu, \dot{\eta}_u + \text{ad}_{\xi}\eta_u - [D_g f(g,u)] g \eta_u \rangle - \langle \mu, \frac{\partial}{\partial u} f(g,u)\rangle \Big] dt \\
&= \langle dC(g(T)), \delta_u g (T)\rangle - \langle \mu, \eta_u\rangle\Big|_{0}^{T} \\ & \qquad \qquad - \int_{0}^{T} \Big[ \langle -\dot{\mu} + \text{ad}^*_{\xi}\mu - g^{*}[D_g f(g,u)]^*\mu, \eta_u \rangle - \langle \mu, \frac{\partial}{\partial u} f(g,u)\rangle \Big]dt,
\end{align*}
where we integrated the $\langle \mu, \dot{\eta}_u\rangle$ term by parts. Now, the first pairing in the integral vanishes if $\mu$ satisfies the adjoint equation. Furthermore, $\eta_u(0) = 0$ since the initial condition for problem \eqref{Parameter Sensitivity Problem} is fixed. Hence, we have
$$ \frac{d}{du}J(g(T)) = \langle dC(g(T)), \delta_u g (T)\rangle - \langle \mu(T), \eta_u(T)\rangle + \int_{0}^{T} \langle \mu, \frac{\partial}{\partial u} f(g,u) \rangle dt. $$
If we choose the terminal condition $\mu(T) = g(T)^* dC(g(T))$, the first two terms on the right hand side cancel, which gives the expression for the desired parameter sensitivity
$$ \frac{d}{du} C(g(T)) = \frac{d}{du} J(g(T)) = \int_{0}^{T} \langle \mu, \frac{\partial}{\partial u}f(g,u)\rangle dt, $$
where $\mu$ is chosen to satisfy the adjoint equation $ -\dot{\mu} + \text{ad}^*_{\xi}\mu - g^{*}[D_g f(g,u)]^*\mu  = 0$ and the terminal condition $\mu(T) = g(T)^* dC(g(T))$.
\end{proof}
\end{prop}

\section{Application to Discrete Variational Integrators}\label{DHVI section}
As an example application of the preceding sections, we develop a discrete counterpart to the continuous Type II variational principle discussed in \Cref{sec:type-ii-vp}. We will derive a first-order Lie group variational integrator for simplicity in presentation, but higher-order schemes can be derived similarly. We will subsequently discuss such integrators in the context of adjoint systems.

Consider the action
$$ s[g,\mu] = \int_0^T \Big( \langle \mu, g^{-1}\cdot \dot{g}\rangle - h(g,\mu) \Big) dt. $$
We will construct discrete Hamiltonian variational integrators for the Lie--Poisson system \eqref{TypeIILPa}-\eqref{TypeIILPd} by discretizing the Type II d'Alembert variational principle (Theorem \ref{TypeIIdAVP Theorem}).

Partition $[0,T]$ into $\cup_{k=0}^{N-1}[t_k,t_{k+1}]$ where 
$$ 0=t_0<t_1<\dots<t_{N-1}<t_N=T, $$
with uniformly spaced intervals $t_{k+1}-t_k = \Delta t = T/N$. To discretize the variational principle, we need a sequence of points $\{g_k \in G\}_{k=0}^{N-1}$ which interpolates a curve $g(t) \in G$. A simple way to do this is to utilize a retraction to relate a curve on $G$ to a curve on $\mathfrak{g}$. Let $\tau$ be a retraction $\tau: \mathfrak{g} \rightarrow G$, which is a $C^2$-diffeomorphism about the origin such that $\tau(0) = e$. Let $d\tau_\xi: \mathfrak{g} \rightarrow \mathfrak{g}$ denote the right-trivialized tangent map of $\tau$ and $d\tau_\xi^{-1}$ its inverse (for a definition, see \cite{BoMa2009}). Using the retraction, we can approximate the velocity $g^{-1}\cdot \dot{g} \in \mathfrak{g}$ by 
\begin{equation}\label{Velocity Approximation}
\xi_{k+1} = \tau^{-1}(g_k^{-1}g_{k+1})/\Delta t.
\end{equation}
This defines the desired interpolated curve $\{g_k\}$ on $G$ through elements $\{\xi_k\}$ on $\mathfrak{g}$ via $g_{k+1} = g_k \tau(\Delta t\xi_{k+1})$. We approximate the action as
$$ s_d[\{g_k\},\{m_k\}] = \sum_{k=0}^{N-1} \Delta t \Big( \langle m_{k+1},\xi_{k+1}\rangle - h(g_k, m_{k+1}) \Big),$$
where again $\{\xi_k\}$ and $\{g_k\}$ are related by \eqref{Velocity Approximation}. Note that, by \eqref{Velocity Approximation}, the variations in $\xi$ are related to the variations of $g$; from \cite{KoMa2011}, this is explicitly given by 
\begin{equation}
\delta\xi_{k+1} = \delta \tau^{-1}(g_k^{-1}g_{k+1})/\Delta t = d\tau^{-1}_{\Delta t\xi_{k+1}}(-g_k^{-1}\delta g_k + \text{Ad}_{\tau(\Delta t\xi_{k+1})}g_{k+1}^{-1}\delta g_{k+1})/\Delta t. \label{gxi variation identity}
\end{equation}
We now derive a variational integrator from a discrete approximation of the Type II d'Alembert variational principle. 

\begin{theorem}[Discrete Type II d'Alembert Variational Principle]\label{DiscreteTypeIIVP Theorem}
The following are equivalent
\begin{enumerate}[label=\normalfont(\roman*)]
\item The discrete Type II d'Alembert variational principle holds
$$ \delta s_d[\{g_k\},\{m_k\}] = \langle (d\tau^{-1}_{-\Delta t \xi_{N}})^*m_N, g_N^{-1} \delta g_N\rangle, $$
subject to variations $\delta g_k, \delta m_k$ satisfying $\delta g_0 = 0, \delta m_N=0$, which correspond to Type II boundary conditions $g_0 = g(0), m_N = m(T)$.
\item The discrete Lie--Poisson equations hold
\begin{subequations}
\begin{align}
(d\tau^{-1}_{\Delta t\xi_{k+1}})^*m_{k+1} - \textup{Ad}^*_{\tau(\Delta t\xi_k)} (d\tau^{-1}_{\Delta t\xi_{k}})^*m_{k} &= - \Delta t g_k^*\cdot D_gh(g_k,m_{k+1}), \label{DiscreteLP a} \\
\xi_{k+1} &= D_\mu h(g_k, m_{k+1}), \label{DiscreteLP b}\\
g_{k+1} &= g_k \tau (\Delta t \xi_{k+1}), \label{DiscreteLP c}
\end{align}
\end{subequations}
with the above boundary conditions.
\end{enumerate}

\begin{proof}
Compute the variation of $s_d$,
\begin{align*}
\delta s_d &= \underbrace{\sum_{k=0}^{N-1} \Delta t \Big[ \langle \delta m_{k+1},\xi_{k+1} \rangle - \langle \delta m_{k+1}, D_\mu h(g_k,m_{k+1}) \rangle \Big] }_{\equiv (a)} \\
& \qquad + \underbrace{\sum_{k=0}^{N-1} \Delta t \Big[ \langle m_{k+1},\delta \xi_{k+1} \rangle - \langle D_gh(g_k,m_{k+1}), \delta g_k \rangle \Big]}_{\equiv (b)}.
\end{align*}
We will simplify the expressions (a) and (b) individually. 

For (a), note that the $k=N-1$ term vanishes since $\delta m_N = 0$. Thus, the sum runs $0$ to $N-2$. We re-index $k \rightarrow k-1$ so that (a) becomes
$$ \text{(a)} = \sum_{k=1}^{N-1} \Delta t \Big[ \langle \delta m_{k},\xi_{k} - D_\mu h(g_{k-1},m_{k}) \rangle \Big]. $$

For (b), we rewrite the variation in $\xi$ in terms of the variation of $g$,
\begin{align*}
\text{(b)} &= \sum_{k=0}^{N-1} \Delta t\Big[ \langle m_{k+1} , d\tau^{-1}_{\Delta t\xi_{k+1}}(-g_k^{-1}\delta g_k + \text{Ad}_{\tau(\Delta t\xi_{k+1})}g_{k+1}^{-1}\delta g_{k+1})/\Delta t\rangle -\langle D_gh(g_k,m_{k+1}), \delta g_k \rangle \Big] \\
&= \sum_{k=0}^{N-1} \Delta t\Big[ \langle -(g_k^{-1})^* (d\tau^{-1}_{\Delta t \xi_{k+1}})^*m_{k+1}/\Delta t - D_gh(g_k,m_{k+1}), \delta g_k\rangle \Big] \\ 
& \qquad \quad + \sum_{k=0}^{N-1} \Delta t \langle (g_{k+1}^{-1})^* \textup{Ad}^*_{\tau(\Delta t\xi_{k+1})} (d\tau^{-1}_{\Delta t{\xi_{k+1}}})^*m_{k+1}/\Delta t, \delta g_{k+1} \rangle.
\end{align*}
In the first sum above, note that the $k=0$ vanishes since $\delta g_0 = 0$. In the second sum above, we explicitly write the $k=N-1$ term and re-index the resulting sum $k \rightarrow k-1$. This gives
\begin{align*}
\text{(b)} &= \sum_{k=1}^{N-1} \Delta t\Big[ \langle -(g_k^{-1})^* (d\tau^{-1}_{\Delta t \xi_{k+1}})^*m_{k+1}/\Delta t - D_gh(g_k,m_{k+1}), \delta g_k\rangle \Big] \\ 
& \qquad \quad + \sum_{k=1}^{N-1} \Delta t \langle (g_{k}^{-1})^* \textup{Ad}^*_{\tau(\Delta t\xi_{k})} (d\tau^{-1}_{\Delta t{\xi_{k}}})^*m_{k}/\Delta t, \delta g_{k} \rangle \\
& \qquad \quad + \Delta t \langle (g_{N}^{-1})^* \textup{Ad}^*_{\tau(\Delta t\xi_{N})} (d\tau^{-1}_{\Delta t{\xi_{N}}})^*m_{N}/\Delta t, \delta g_{N} \rangle.
\end{align*}
Note that, since $\text{Ad}^*_{\tau(\Delta t \xi_N)}(d\tau^{-1}_{\Delta t \xi_N})^* = (d\tau^{-1}_{-\Delta t \xi_N})^*$ \cite{BoMa2009}, the last term equals 
$$\langle (d\tau^{-1}_{-\Delta t \xi_{N}})^*m_N, g_N^{-1} \delta g_N\rangle,$$ which is precisely the virtual work term in the discrete Type II d'Alembert variational principle. Putting everything together, we have
\begin{align*}
 \delta s_d&[\{g_k\},\{m_k\}] - \langle (d\tau^{-1}_{-\Delta t \xi_{N}})^*m_N, g_N^{-1} \delta g_N\rangle \\
&= \sum_{k=1}^{N-1} \Delta t \Big[ \langle \delta m_{k},\xi_{k} - D_\mu h(g_{k-1},m_{k}) \rangle \\
& \qquad \qquad + \langle -(g_k^{-1})^* (d\tau^{-1}_{\Delta t \xi_{k+1}})^*m_{k+1}/\Delta t - D_gh(g_k,m_{k+1}) \\ 
& \qquad \qquad + (g_{k}^{-1})^* \textup{Ad}^*_{\tau(\Delta t\xi_{k})} (d\tau^{-1}_{\Delta t{\xi_{k}}})^*m_{k}/\Delta t, \delta g_k\rangle\Big].
\end{align*}
Clearly, if the discrete Lie--Poisson equations hold, then the above vanishes, i.e., the discrete Type II d'Alembert variational principle holds, $\delta s_d[\{g_k\},\{m_k\}] = \langle (d\tau^{-1}_{-\Delta t \xi_{N}})^*m_N, g_N^{-1} \delta g_N\rangle.$ Conversely, if the discrete Type II d'Alembert variational principle holds, the above vanishes, which gives 
\begin{align*}
0 &= \sum_{k=1}^{N-1} \Delta t \Big[ \langle \delta m_{k},\xi_{k} - D_\mu h(g_{k-1},m_{k}) \rangle \\
& \qquad \qquad + \langle -(g_k^{-1})^* (d\tau^{-1}_{\Delta t \xi_{k+1}})^*m_{k+1}/\Delta t - D_gh(g_k,m_{k+1})  \\
& \qquad \qquad + (g_{k}^{-1})^* \textup{Ad}^*_{\tau(\Delta t\xi_{k})} (d\tau^{-1}_{\Delta t{\xi_{k}}})^*m_{k}/\Delta t, \delta g_k\rangle\Big].
\end{align*}
Since the variations $\delta m_k$ and $\delta g_k$ are arbitrary and independent for $k=1,\dots,N-1$, this gives the discrete Lie--Poisson equations \eqref{DiscreteLP a}-\eqref{DiscreteLP c}.
\end{proof}

\end{theorem}

This is an example of a first-order Lie group variational integrator (see, for example, \cite{BoMa2009, MaPeSh1999, MaRo2010, KoMa2011}), derived using a Type II variational principle and particularly, without making use of a Lagrangian description, so it applies to generic (possibly degenerate) Hamiltonian systems on $T^*G$.

\begin{remark}
In the case that the Hamiltonian is left-invariant, it is straightforward to verify that the above system \eqref{DiscreteLP a}-\eqref{DiscreteLP c} reduces to
\begin{subequations}\label{DiscreteLP reduced}
\begin{align}
(d\tau^{-1}_{\Delta t\xi_{k+1}})^*m_{k+1} - \textup{Ad}^*_{\tau(\Delta t\xi_k)} (d\tau^{-1}_{\Delta t\xi_{k}})^*m_{k} &= 0, \label{DiscreteLP reduced a} \\
\xi_{k+1} &= D_\mu \tilde{H}(m_{k+1}), \label{DiscreteLP reduced b}\\
g_{k+1} &= g_k \tau (\Delta t \xi_{k+1}). \label{DiscreteLP reduced c}
\end{align}
\end{subequations}
\end{remark}

\begin{remark}
It is interesting to note that the virtual work term arising at the terminal point in the discrete variational principle,
$$ \langle (d\tau^{-1}_{-\Delta t\xi_N})^*m_N, g_N^{-1}\delta g_N\rangle $$
is different than what one might expect from the continuous variational principle, $\langle m_N, g_N^{-1} \delta g_N\rangle$. This is due to the fact that the retraction relates the dynamics on $G$ to dynamics on $\mathfrak{g}$, and so the pairing $\langle \mu, g^{-1}\delta g\rangle$ compared to the pairing $\langle m_k,g_k^{-1} \delta g_k\rangle$ should not be identified, but rather, are related by a coordinate change; namely, $\mu = (d\tau^{-1}_{-\Delta t\xi}) m$. In fact, this coordinate change is given by the cotangent lift of $\tau^{-1}$, which is precisely $(d\tau^{-1}_{-\Delta t\xi_N})^*$. As we will see below, this also induces a coordinate change in the expression for the symplectic form, which is the exterior derivative of the one-form corresponding to the above boundary term; the expression for the one-form and its exterior derivative is also derived in \cite{BoMa2009} through a discrete Hamilton--Pontryagin principle.
\end{remark}

\textbf{Discrete Conservation Properties.} It is straightforward (though tedious) to verify that the above method satisfies a discrete symplectic conservation law and a discrete Noether theorem. For brevity, we will simply state the results here; the proofs follow from a direct calculation.

From the boundary term arising from the variation of the discrete action $s_d$, we see that the discrete canonical form has the expression
\begin{equation}\label{Discrete Canonical Form}
\Theta_k = \langle (d\tau^{-1}_{-\Delta t \xi_k})^*m_k, g_k^{-1}dg_k\rangle,
\end{equation}
whose action on a vector $\delta_k = \delta m_k \partial/\partial m_k + \delta g_k \partial/\partial g_k$ is given by
$$ \Theta_k \cdot \delta_k = \langle (d\tau^{-1}_{-\Delta t \xi_k})^*m_k, g_k^{-1}\delta g_k\rangle. $$
Then, the corresponding discrete symplectic form $\Omega_k \equiv d\Theta_k$ has the expression
\begin{equation}\label{Discrete Symplectic Form}
\Omega_k = (d\tau^{-1}_{-\Delta t\xi_k})^* dm_k \wedge g_k^{-1}dg_k.
\end{equation}
Its action on vectors $\delta^i_k = \delta m^i_k \partial/\partial m_k + \delta g^i_k \partial/\partial g_k$ is given by
$$ \Omega_k\cdot (\delta^1_k, \delta^2_k) = \langle (d\tau^{-1}_{-\Delta t \xi_k})^* \delta m^1_k, g_k^{-1}\delta g^2_k\rangle - \langle (d\tau^{-1}_{-\Delta t \xi_k})^* \delta m^2_k, g_k^{-1}\delta g^1_k\rangle. $$
Symplecticity of the integrator \eqref{DiscreteLP a}-\eqref{DiscreteLP c} is the statement that $\Omega_{k+1} = \Omega_k$ when the discrete Lie--Poisson equations \eqref{DiscreteLP a}-\eqref{DiscreteLP c} hold, where the symplectic forms are evaluated on first variations of the discrete Lie--Poisson equations, i.e., variations whose flow preserves solutions of the discrete Lie--Poisson equations. Equivalently, such first variations are those which preserve \eqref{DiscreteLP a}-\eqref{DiscreteLP c} to linear order. 

\begin{theorem}\label{Discrete Symplecticity Theorem}
The integrator \eqref{DiscreteLP a}-\eqref{DiscreteLP c} is symplectic, i.e., the symplectic form is preserved,
$$\Omega_{k+1} = \Omega_k,$$
subject to first variations of the discrete Lie--Poisson equations.
\end{theorem}

The integrator \eqref{DiscreteLP a}-\eqref{DiscreteLP c} also satisfies a discrete Noether's theorem. Let $\{g_k^\epsilon, m_k^\epsilon\}$ be a one-parameter family of discrete time curves with $g_k^0 = g_k$ and $m_k^0 = m_k$. Denote by
\begin{align*}
\delta g_k &= \frac{d}{d\epsilon}g_k^\epsilon \Big|_{\epsilon = 0}, \\
\delta m_k &= \frac{d}{d\epsilon}m_k^\epsilon \Big|_{\epsilon = 0},
\end{align*}
the variations associated with the one-parameter family of discrete time curves. Furthermore, let $s_k = \langle m_{k+1},\xi_{k+1}\rangle - h(g_k, m_{k+1})$ denote the $k^{th}$ discrete action density. Then, we have the following momentum preservation property of \eqref{DiscreteLP a}-\eqref{DiscreteLP c}.

\begin{theorem}[Discrete Noether's Theorem]\label{Discrete Momentum Conservation Theorem}
Suppose that \eqref{DiscreteLP a}-\eqref{DiscreteLP c} hold and furthermore, suppose that the $k^{th}$ discrete action density is invariant under the above variations,
$$ \delta s_k = 0.$$
Then, for any time indices $I < J$, 
\begin{equation}\label{Discrete Momentum Conservation Equation}
\Theta_I \cdot \delta g_I = \Theta_J \cdot \delta g_J,
\end{equation}
where $\Theta_k$ is the discrete canonical form \eqref{Discrete Canonical Form}.
\end{theorem}

As an application of Theorem \ref{Discrete Momentum Conservation Theorem}, we will re-derive the discrete reduced Lie--Poisson equation \eqref{DiscreteLP reduced}, interpreted as momentum conservation associated with left-invariance symmetry. Let $H$ be a left-invariant Hamiltonian, let $X$ be a right-invariant vector field on $G$ with $X(e) = \chi \in \mathfrak{g}$, and let $\varphi_\epsilon$ denote the time-$\epsilon$ flow of $X$. We choose $X$ to be a right-invariant vector field, since its flow is given by left translations
$$ \varphi_\epsilon(g) = \exp(\epsilon \chi) g, $$
where $\exp: \mathfrak{g} \rightarrow G$ denotes the exponential map. We define a one-parameter family of discrete time curves $\{g_k^\epsilon, m_k^\epsilon\}$ as 
\begin{align*}
g_k^\epsilon &= \varphi_\epsilon (g_k) = \exp(\epsilon \chi)g, \\
m_k^\epsilon &= m_k,
\end{align*}
i.e., the one-parameter family of discrete time curves is defined by flowing $g_k$ by $\varphi_\epsilon$, whereas $m_k^\epsilon$ remains constant with $\epsilon$. To see why we defined $m_k^\epsilon$ this way, since $(g_k, p_k)$ transforms under the cotangent lift of left-multiplication by $x$ as $(g_k,p_k) \mapsto (xg_k, x^{*-1}p_k)$, $m_k$ transforms as
$$ m_k = g_k^* p_k \mapsto (xg_k)^* x^{*-1}p_k = g_k^* x^* x^{*-1} p_k = g_k^*p_k = m_k, $$
i.e., $m_k$ is invariant under this transformation; thus, we define $m_k^\epsilon$ to be constant in $\epsilon$. Additionally, observe that the variations associated with this one-parameter family of discrete time curves can be expressed as
\begin{align*}
\delta g_k &= \frac{d}{d\epsilon}\Big|_0 g_k^\epsilon = \frac{d}{d\epsilon}\Big|_0 \varphi_\epsilon(g_k) = X(g_k), \\
\delta m_k &= \frac{d}{d\epsilon}\Big|_0 m_k^\epsilon = \frac{d}{d\epsilon}\Big|_0 m_k = 0.
\end{align*}

Now, we will verify the assumption of Theorem \ref{Discrete Momentum Conservation Theorem}. The $k^{th}$ discrete action density is
$$ s_k = \langle m_{k+1},\xi_{k+1}\rangle - h(g_k,m_{k+1}) = \langle m_{k+1}, \xi_{k+1}\rangle - \tilde{H}(m_{k+1}), $$
where in the second equality, we used that $h(g_k,m_{k+1}) = \tilde{H}(m_{k+1})$ for a left-invariant Hamiltonian. As stated above, $m_{j}^\epsilon = m_j$ is invariant under this transformation. Furthermore, since $\xi_j = \tau^{-1}(g_k^{-1}g_{k+1})/\Delta t$, the corresponding transformation for $\xi_j$ is given by 
\begin{align*}
\xi_j^\epsilon &= \tau^{-1}( (g_k^\epsilon)^{-1}g_{k+1}^{\epsilon} )/\Delta = \tau^{-1}( (\exp(\epsilon \chi)g_k)^{-1}\exp(\epsilon \chi)g_{k+1} )/\Delta t \\
&= \tau^{-1}( g_k^{-1} \exp(\epsilon \chi)^{-1}\exp(\epsilon \chi)g_{k+1} )/\Delta = \tau^{-1}(g_k^{-1}g_{k+1})/\Delta t = \xi_j,
\end{align*}
i.e., $\xi_j$ is also invariant under this transformation. Hence, $s_k$ is invariant under the above variation, so Theorem \ref{Discrete Momentum Conservation Theorem} applies. We thus have $\Theta_{k+1}\cdot \delta g_{k+1} = \Theta_k \cdot \delta g_k$, i.e.,
$$ \langle (d\tau^{-1}_{-\Delta t \xi_{k+1}})^*m_{k+1}, g_{k+1}^{-1}\delta g_{k+1}\rangle = \langle (d\tau^{-1}_{-\Delta t \xi_{k}})^*m_{k}, g_{k}^{-1}\delta g_{k}\rangle .$$
It is straightforward to verify that this is equivalent to \eqref{DiscreteLP reduced a}.

\subsection{Discrete Adjoint Sensitivity Analysis on Lie Groups}\label{Adjoint Sensitivity Analysis Numerical Section}
In this section, we apply the Type II variational integrators developed in Section \ref{DHVI section} to the particular case of adjoint systems. We will show explicitly that these integrators preserve the adjoint-variational quadratic conservation law which gives rise to the adjoint sensitivity method discussed in \Cref{Adjoint Sensitivity Analysis Section}, and thus, these methods are geometric structure-preserving methods for adjoint sensitivity analysis on Lie groups. For discussions on the application of variational integrators to geometric optimal control, see for example \cite{Sig2024, Col2016}.

Consider the variational integrators that we derived in Section \ref{DHVI section}, applied to the adjoint system \eqref{LP Adjoint a}-\eqref{LP Adjoint b}. Substituting $h(g,\mu) = \langle\mu,f(g)\rangle$ into the discrete Lie--Poisson equations \eqref{DiscreteLP a}-\eqref{DiscreteLP c}, we have the discrete Lie--Poisson adjoint equations
\begin{subequations}
\begin{align}
(d\tau^{-1}_{\Delta t \xi_{k+1}})^*m_{k+1} - \text{Ad}^*_{\tau(\Delta t \xi_k)}(d\tau^{-1}_{\Delta t \xi_k})^*m_k &= - \Delta t g_k^* [Df(g_k)]^*m_{k+1}, \label{Discrete LP Adjoint a} \\
\xi_{k+1} &= f(g_k) \label{Discrete LP Adjoint b} \\
g_{k+1} &= g_k \tau(\Delta t\xi_{k+1}) = g_k \tau(\Delta t f(g_k)). \label{Discrete LP Adjoint c}
\end{align}
\end{subequations}

In order to derive a discrete analogue of the adjoint conservation law, we consider the discrete variational equation, which is a discretization of the continuous variational equation \eqref{VariationalEquation Left-Trivial b}. To derive the discrete variational equation, note that as mentioned in Section \ref{DHVI section}, the variation of equation \eqref{Discrete LP Adjoint c} can be expressed as
$$ \delta \xi_{k+1} = d\tau^{-1}_{\Delta t \xi_{k+1}}(-g_k^{-1}\delta g_k + \text{Ad}_{\tau(\Delta t \xi_{k+1})}g_{k+1}^{-1}\delta g_{k+1})/\Delta t. $$
Furthermore, by taking the variation of equation \eqref{Discrete LP Adjoint b}, we have
$$ \delta \xi_{k+1} = Df(g_k)\delta g_k. $$
Combining these two equations yields
$$ Df(g_k)\delta g_k = d\tau^{-1}_{\Delta t \xi_{k+1}}(-g_k^{-1}\delta g_k + \text{Ad}_{\tau(\Delta t \xi_{k+1})}g_{k+1}^{-1}\delta g_{k+1})/\Delta t. $$
Defining the left-trivialized variation $\eta_k = g_k^{-1}\delta g_k$, the above can be expressed as
\begin{equation}\label{DiscreteVariationalEquation}
\Delta t d\tau_{\Delta t\xi_{k+1}}Df(g_k)g_k\eta_k = - \eta_k + \text{Ad}_{\tau(\Delta t \xi_{k+1})}\eta_{k+1},
\end{equation}
which we refer to as the discrete variational equation.

\begin{theorem}\label{Discrete Quadratic Conservation Law Theorem}
The discrete Lie--Poisson adjoint equations and the discrete variational equation satisfy the following quadratic conservation law,
\begin{equation}\label{DiscreteAdjointConservationLaw}
\langle (d\tau^{-1}_{-\Delta t \xi_{k+1}})^*m_{k+1}, \eta_{k+1}\rangle = \langle (d\tau^{-1}_{-\Delta t \xi_k})^*m_k,\eta_k\rangle.
\end{equation}
\begin{proof}
We use the identity $\text{Ad}^*_{\tau(\Delta t\xi_j)}(d\tau^{-1}_{\Delta t \xi_j})^* = (d\tau^{-1}_{-\Delta t \xi_j})^*$. Starting from the left hand side of equation \eqref{DiscreteAdjointConservationLaw}, we compute
\begin{align*}
\langle (d\tau^{-1}_{-\Delta t \xi_{k+1}})^*m_{k+1}, \eta_{k+1}\rangle &=  \langle \text{Ad}^*_{\tau(\Delta t\xi_{k+1})}(d\tau^{-1}_{\Delta t \xi_{k+1}})^*m_{k+1}, \eta_{k+1}\rangle \\
&= \langle (d\tau^{-1}_{\Delta t \xi_{k+1}})^*m_{k+1}, \text{Ad}_{\tau(\Delta t\xi_{k+1})}\eta_{k+1}\rangle.
\end{align*}
Substituting \eqref{Discrete LP Adjoint a} and \eqref{DiscreteVariationalEquation} yields
\begin{align*}
\langle (d\tau^{-1}_{-\Delta t \xi_{k+1}}&)^* m_{k+1}, \eta_{k+1}\rangle\\
& = \langle (d\tau^{-1}_{\Delta t \xi_{k+1}})^*m_{k+1}, \text{Ad}_{\tau(\Delta t\xi_{k+1})}\eta_{k+1}\rangle \\
&= \langle (d\tau^{-1}_{-\Delta t \xi_k})^*m_k - \Delta t g_k^* [Df(g_k)]^*m_{k+1}, \eta_k + \Delta t d\tau_{\Delta t\xi_{k+1}}Df(g_k)g_k\eta_k \rangle \\
&= \langle (d\tau^{-1}_{-\Delta t \xi_k})^*m_k,\eta_k\rangle + \Delta t \langle (d\tau^{-1}_{-\Delta t \xi_k})^*m_k, d\tau_{\Delta t \xi_{k+1}}Df(g_k)g_k\eta_k\rangle \\
& \quad - \Delta t \langle g_k^* [Df(g_k)]^*m_{k+1},\eta_k\rangle - \Delta t^2 \langle g_k^* [Df(g_k)]^*m_{k+1}, d\tau_{\Delta t \xi_{k+1}}Df(g_k)g_k\eta_k\rangle.
\end{align*}
Substitute \eqref{Discrete LP Adjoint a} into the second term above,
\begin{align*}
\langle (d&\tau^{-1}_{-\Delta t \xi_{k+1}})^* m_{k+1}, \eta_{k+1}\rangle \\
&= \langle (d\tau^{-1}_{-\Delta t \xi_k})^*m_k,\eta_k\rangle + \Delta t \langle (d\tau^{-1}_{\Delta t \xi_{k+1}})^*m_{k+1} + \Delta t g_k^* [Df(g_k)]^*m_{k+1}, d\tau_{\Delta t \xi_{k+1}}Df(g_k)g_k\eta_k\rangle \\
& \quad  - \Delta t \langle g_k^* [Df(g_k)]^*m_{k+1},\eta_k\rangle - \Delta t^2 \langle g_k^* [Df(g_k)]^*m_{k+1}, d\tau_{\Delta t \xi_{k+1}}Df(g_k)g_k\eta_k\rangle \\
&= \langle (d\tau^{-1}_{-\Delta t \xi_k})^*m_k,\eta_k\rangle +\Delta t \langle (d\tau^{-1}_{\Delta t \xi_{k+1}})^*m_{k+1}, d\tau_{\Delta t \xi_{k+1}}Df(g_k)g_k\eta_k\rangle \\
& \quad + \Delta t^2 \langle g_k^* \cancel{[Df(g_k)]^*m_{k+1}, d\tau_{\Delta t \xi_{k+1}} }Df(g_k)g_k\eta_k\rangle \\
& \quad  - \Delta t \langle g_k^* [Df(g_k)]^*m_{k+1},\eta_k\rangle - \Delta t^2 \langle g_k^* \cancel{[Df(g_k)]^*m_{k+1}, d\tau_{\Delta t \xi_{k+1}} }Df(g_k)g_k\eta_k\rangle \\
&= \langle (d\tau^{-1}_{-\Delta t \xi_k})^*m_k,\eta_k\rangle +\Delta t \langle m_{k+1}, (d\tau^{-1}_{\Delta t \xi_{k+1}}) d\tau_{\Delta t \xi_{k+1}}Df(g_k)g_k\eta_k\rangle \\
& \quad  - \Delta t \langle g_k^* [Df(g_k)]^*m_{k+1},\eta_k\rangle \\
&= \langle (d\tau^{-1}_{-\Delta t \xi_k})^*m_k,\eta_k\rangle +\Delta t \langle m_{k+1}, Df(g_k)g_k\eta_k\rangle - \Delta t \langle g_k^* [Df(g_k)]^*m_{k+1},\eta_k\rangle \\
&= \langle (d\tau^{-1}_{-\Delta t \xi_k})^*m_k,\eta_k\rangle. \qedhere
\end{align*}
\end{proof}
\end{theorem}

We now consider discrete counterparts of the optimization problems discussed in \Cref{Adjoint Sensitivity Analysis Section}.

\textbf{Initial Condition Sensitivity.} We begin with problem \eqref{Adjoint Sensitivity Problem} in the discrete setting. Analogous to the continuous case discussed in \eqref{Adjoint Sensitivity Analysis Section}, the initial condition sensitivity can be obtained from the associated quadratic conservation law for the adjoint-variational system discussed above. More specifically, to obtain the initial condition sensitivity, recall the quadratic conservation law
$$ \langle (d\tau^{-1}_{-\Delta t \xi_{N}})^*m_{N}, g_N^{-1}\delta g_N\rangle = \langle (d\tau^{-1}_{-\Delta t \xi_0})^*m_0,g_0^{-1} \delta g_0\rangle. $$
We set $(d\tau^{-1}_{-\Delta t \xi_{N}})^*m_{N} = d_L C(g_N)$, where $d_L$ denotes the left-trivialized derivative, $d_L C(g_N) \equiv g_N^{*} dC(g_N)$. This gives $m_N = (d\tau_{-\Delta t \xi_{N}})^* d_LC(g_N)$. Subsequently, evolve the momenta backward in time using \eqref{Discrete LP Adjoint a} to obtain $m_0$. Finally, the left-trivialized derivative of $C(g_N)$ with respect to $g_0$ is given by $(d\tau^{-1}_{-\Delta t \xi_0})^*m_0$, which follows from a similar calculation to the continuous case. This is summarized in \Cref{Discrete Gradient Algorithm}.

\begin{algorithm}
\caption{Left-Trivialized Initial Condition Sensitivity}
\label{Discrete Gradient Algorithm}
\begin{algorithmic}
\State \textbf{Input:}	$g_{\text{init}}$
\State \textbf{Initialize:} $g_0 \gets g_{\text{init}}$
\State \textbf{Output:} Left-Trivialized Derivative of $C(g_N)$ with respect to $g_0$
\For{k=1,...,N}
\State $g_k \gets g_{k-1} \tau(\Delta t f(g_{k-1}))$
\State $\xi_k \gets \tau^{-1}(g^{-1}_{k-1}g_k)/\Delta t$
\EndFor
\State $m_N \gets (d\tau_{-\Delta t f(g_N)})^* d_LC(g_N)$
\For{k=N-1,\dots,0}
\State \text{Solve \eqref{Discrete LP Adjoint a} for} $m_k$
\EndFor
\State \textbf{Return } $(d\tau^{-1}_{-\Delta t f(g_0)})^*m_0$
\end{algorithmic}
\end{algorithm}

This can be combined with a line-search algorithm to solve the optimization problem \eqref{Adjoint Sensitivity Problem}. More precisely, fixing an inner product on $\mathfrak{g}$, such as the Frobenius inner product
$$ (A,B)_F \equiv \text{Tr}(A^*B), $$
we can identify $\mathfrak{g}^*$ with $\mathfrak{g}$ and hence, identify the output of Algorithm \ref{Discrete Gradient Algorithm} with the discrete left-trivialized gradient of $C(g_N)$ with respect to $g_0$, $\widetilde{\nabla}_{g_0}C(g_N)$, which is an element of $\mathfrak{g}$. With this identification, the initial condition can be updated as $g_0 \leftarrow g_0 \tau(- \gamma \widetilde{\nabla}_{g_0}C(g_N))$, for some line-search step size $\gamma$. Note that this approach is intrinsic: at any iteration in the line-search algorithm, the iterate $g_0$ is valued in $G$, to numerical precision. Furthermore, while this is also true of projection-based optimization algorithms, such methods generally no longer preserve the adjoint-variational quadratic conservation law and hence, may yield search directions that are not guaranteed to be descent directions. More precisely, by satisfying the discrete quadratic conservation law, the sensitivity produced from \Cref{Discrete Gradient Algorithm} produces the exact sensitivity for the discrete cost function $C(g_N)$ \cite{TrLe2022adj}.

\textbf{Parameter Sensitivity.} We now consider the discrete setting for the optimal control problem problem \eqref{Parameter Sensitivity Problem} discussed in \Cref{Adjoint Sensitivity Analysis Section}. In analogy with the continuous case, we define the parameter-dependent left-trivialized discrete action
$$s_d[\{g_k\},\{m_k\};u] = \Delta t \sum_{k=0}^{N-1} \Big( \langle m_{k+1},\xi_{k+1}\rangle - h(g_k,m_{k+1}l; u) \Big), $$
and form the discrete augmented cost function
\begin{align*}
J_d &\equiv C(g_N) - s_d[\{g_k\},\{m_k\};u] = C(g_N) - \Delta t \sum_{j=0}^{N-1} \langle m_{j+1}, \xi_{j+1} - f(g_j,u)\rangle. 
\end{align*}
We then have an analogous result to determine the parameter sensitivity by computing the derivative $dJ_d/du$.
\begin{prop}
The discrete parameter sensitivity is given by
\begin{equation}\frac{d}{du} C(g_N) = \Delta t \sum_{j=0}^{N-1} \left\langle m_j, \frac{\partial}{\partial u}f(g_j,u)\right\rangle,
\end{equation}
where $m_j$ is chosen to satisfy the discrete Lie--Poisson adjoint equation  \eqref{Discrete LP Adjoint a} with $f(g_k)$ replaced by $f(g_k,u)$, i.e.,
\begin{equation}\label{Discrete LP adjoint a parameter}
    (d\tau^{-1}_{\Delta t \xi_{k+1}})^*m_{k+1} - \text{Ad}^*_{\tau(\Delta t \xi_k)}(d\tau^{-1}_{\Delta t \xi_k})^*m_k = - \Delta t g_k^* [Df(g_k,u)]^*m_{k+1}, 
\end{equation}
and the terminal condition $m_N =(d\tau_{-\Delta t \xi_N})^* d_LC(g_N)$.
\begin{proof}
Analogous to the continuous setting, we have
$$ \frac{d}{du} C(g_N) = \frac{d}{du} J_d. $$
Now, we calculate $dJ_d/du$ explicitly,
\begin{align*}
\frac{d}{du}J_d &= \langle dC(g_N), \delta_u g_N\rangle - \Delta t \sum_{j=0}^{N-1} \left\langle m_{j+1}, \delta_u \xi_{j+1} - D_g f(g_j,u) \delta_u g_j - \frac{\partial f}{\partial u}(g_j,u) \right\rangle.
\end{align*}
Using the identity $\delta_u \xi_{j+1} = d\tau^{-1}_{\Delta t \xi_{j+1}}(-g_k^{-1} \delta_u g_k + \text{Ad}_{\tau(\Delta t \xi_{j+1})}g_{j+1}^{-1}\delta_u g_{j+1})/\Delta t$, the above can be expressed as
\begin{align*}
\frac{d}{du}J_d &= \langle dC(g_N), \delta_u g_N\rangle - \Delta t \sum_{j=0}^{N-1} \Big[ \Delta t^{-1} \langle (d\tau^{-1}_{\Delta t \xi_{j+1}})^*m_{j+1}, -g_j^{-1}\delta_u g_j\rangle \\
 &\quad +\Delta t^{-1} \langle\text{Ad}^*_{\tau(\Delta t \xi_{j+1})}(d\tau^{-1}_{\Delta t \xi_{j+1}})^*m_{j+1}, g_{j+1}^{-1}\delta_u g_{j+1} \rangle \\
& \quad - \langle g_j^{*}[D_gf(g_j,u)]^*m_{j+1}, g_j^{-1}\delta_u g_j\rangle - \left\langle m_{j+1},\frac{\partial f}{\partial u}(g_j,u)\right\rangle\Big].
\end{align*}
Now, we reindex $j \rightarrow j-1$ the second pairing inside the square brackets above; the sum for this term now runs from $1$ to $N$. However, we explicitly write out the $j=N$ term and note that we can include the $j=0$ term in the sum since $\delta_u g_0 = 0$, as the initial condition $g_0$ is fixed under the variation. Hence, we have
\begin{align*}
\frac{d}{du}J_d &= \langle dC(g_N), \delta_u g_N\rangle - \langle (d\tau^{-1}_{-\Delta t \xi_N})^*m_N, g_N^{-1}\delta_u g_N \rangle \\
& \quad -\Delta t \sum_{j=0}^{N-1} \Big[ \Delta t^{-1} \langle (d\tau^{-1}_{\Delta t \xi_{j+1}})^*m_{j+1}, -g_j^{-1}\delta_u g_j\rangle \\ 
& \quad +\Delta t^{-1} \langle\text{Ad}^*_{\tau(\Delta t \xi_{j})}(d\tau^{-1}_{\Delta t \xi_{j}})^*m_{j}, g_{j}^{-1}\delta_u g_{j} \rangle 
\\
& \quad - \langle g_j^{*}[D_gf(g_j,u)]^*m_{j+1}, g_j^{-1}\delta_u g_j\rangle - \left\langle m_{j+1},\frac{\partial f}{\partial u}(g_j,u)\right\rangle\Big].
\end{align*}
The first two terms above vanish if we set the terminal condition $(d\tau^{-1}_{-\Delta t \xi_N})^*m_N = g_N^{*}dC(g_N)$, i.e., $m_N =(d\tau_{-\Delta t \xi_N})^* d_LC(g_N)$. Furthermore, the first three terms in the square brackets vanish if $m_j$ satisfies the discrete Lie--Poisson adjoint equation \eqref{Discrete LP adjoint a parameter}. Hence, we have the parameter sensitivity
\[ \frac{d}{du} C(g_N) = \frac{d}{du}J_d = \Delta t \sum_{j=0}^{N-1} \left\langle m_{j+1}, \frac{\partial f}{\partial u}(g_j,u)\right\rangle. \qedhere \]
\end{proof}
\end{prop}

The parameter sensitivity is summarized in the following algorithm.
\begin{algorithm}[H]
\caption{Parameter Sensitivity}
\label{Discrete Parameter Gradient Algorithm}
\begin{algorithmic}
\State \textbf{Input:}	$g_{\text{init}}$, $u_{\text{init}}$
\State \textbf{Initialize:} $g_0 \gets g_{\text{init}}$, $u \gets u_{\text{init}}$
\State \textbf{Output:} Derivative of $C(g_N)$ with respect to $u$
\For{k=1,...,N}
\State $g_k \gets g_{k-1} \tau(\Delta t f(g_{k-1},u))$
\State $\xi_k \gets \tau^{-1}(g^{-1}_{k-1}g_k)/\Delta t$
\EndFor
\State $m_N \gets (d\tau_{-\Delta t f(g_N,u)})^* d_LC(g_N)$
\For{k=N-1,...,0}
\State \text{Solve \eqref{Discrete LP adjoint a parameter} for} $m_k$
\EndFor
\State \textbf{Return } $\Delta t \sum_{j=0}^{N-1} \langle m_{j+1}, \frac{\partial}{\partial u}f(g_j,u)\rangle$
\end{algorithmic}
\end{algorithm}

This can be combined with a line-search algorithm to solve the optimization problem \eqref{Parameter Sensitivity Problem}. Note that $U$ could be a vector space, in which case a standard line-search algorithm could be used, or $U$ could be a manifold, in which case a line-search algorithm on manifolds could be used (see, for example, \cite{AbMaSe2008}).

\section{Conclusion}
In this paper, we developed a global Type II variational principle on Lie groups and discussed adjoint systems on the trivialization of the cotangent bundle of Lie groups. As an application, we discussed how to use this variational principle to construct discrete variational integrators and subsequently how this can be used to perform structure-preserving adjoint sensitivity analysis on Lie groups, allowing one to exactly compute sensitivities in optimization problems subject to the dynamics of an ODE on $G$.

It would be interesting to synthesize the ideas in \cite{TrLe2022adj}, where we discussed properties of adjoint systems for evolutionary PDEs on Banach spaces, with the ideas presented here, to develop Hamiltonian variational principles and integrators for PDEs where the solutions are valued in Lie groups, algebras, or more generally, solutions which are stationary sections over principal and fibre bundles associated with a structure group $G$, such as gauge field theories (see, for example, \cite{Na2003, Ha2017}). It would be particularly interesting to extend the Type II multisymplectic Hamiltonian variational integrators developed in \cite{TrLe2022} to apply to the setting of Lie group-valued fields, in order to investigate the role of multisymplectic integrators for adjoint sensitivity analysis in both space and time.

Another natural research direction would be to explore the applications of geometric structure-preserving adjoint sensitivity analysis on Lie groups. One such application is the training of neural networks via backpropagation. In particular, if a neural network is viewed as a discretization of a \textit{neural ODE} \cite{ChRuBeDu2018}, then backpropagation can be viewed as a discretization of the corresponding adjoint equation \cite{matsubara2021symplectic}. As is discussed in \cite{matsubara2021symplectic}, utilizing symplectic methods to perform backpropagation leads to efficient methods for training neural networks. It would be interesting to explore symplectic backpropagation of neural networks where the neural ODE arises from a group-equivariant neural network ~\cite{CoWe2016, KoTr2018} where a Lie group symmetry is a fundamental feature of the neural network. In particular, the reduction theory for adjoint systems on Lie groups that was developed in this paper would be relevant.

\section*{Funding}

BKT was supported by the Marc Kac Postdoctoral Fellowship at the Center for Nonlinear Studies at Los Alamos National Laboratory. ML was supported in part by NSF under grants DMS-1345013, CCF-2112665, DMS-2307801, by AFOSR under grant FA9550-23-1-0279. Los Alamos National Laboratory report LA-UR-23-32526.


\bibliographystyle{plainnat}
\bibliography{hvi_lie.bib}

\begin{thebibliography}{42}
\providecommand{\natexlab}[1]{#1}
\providecommand{\url}[1]{\texttt{#1}}
\expandafter\ifx\csname urlstyle\endcsname\relax
  \providecommand{\doi}[1]{doi: #1}\else
  \providecommand{\doi}{doi: \begingroup \urlstyle{rm}\Url}\fi

\bibitem[A.~Lasota(1967)]{LaOp1967}
Z.~Opial A.~Lasota.
\newblock On the existence and uniqueness of solutions of a boundary value problem for an ordinary second-order differential equation.
\newblock \emph{Colloquium Mathematicae}, 18\penalty0 (1):\penalty0 1--5, 1967.

\bibitem[Absil et~al.(2008)Absil, Mahony, and Sepulchre]{AbMaSe2008}
P.-A. Absil, R.~Mahony, and R.~Sepulchre.
\newblock \emph{Optimization Algorithms on Matrix Manifolds}.
\newblock Princeton University Press, 2008.

\bibitem[Bailey and Shampine(1969)]{BaSh1969}
Paul~B Bailey and L.F Shampine.
\newblock Existence from uniqueness for two point boundary value problems.
\newblock \emph{Journal of Mathematical Analysis and Applications}, 25\penalty0 (3):\penalty0 569--574, 1969.
\newblock ISSN 0022-247X.
\newblock \doi{10.1016/0022-247X(69)90256-X}.

\bibitem[Benning et~al.(2019)Benning, Celledoni, Ehrhardt, Owren, and Sch\"onlieb]{DeCeEhOwSc2019}
M.~Benning, E.~Celledoni, {M.~J.} Ehrhardt, B.~Owren, and C.-B. Sch\"onlieb.
\newblock Deep learning as optimal control problems: Models and numerical methods.
\newblock \emph{J. Comput. Dyn.}, 6\penalty0 (2):\penalty0 171--198, 2019.

\bibitem[Bobenko and Suris(1999)]{BoSu1999}
A.~I. Bobenko and Yu.~B. Suris.
\newblock Discrete time {L}agrangian mechanics on {L}ie groups, with an application to the {L}agrange top.
\newblock \emph{Commun. Math. Phys.}, 204:\penalty0 147--188, 1999.

\bibitem[Bou-Rabee and Marsden(2009)]{BoMa2009}
N.~Bou-Rabee and J.~E. Marsden.
\newblock {H}amilton--{P}ontryagin integrators on {L}ie groups part {I}: Introduction and structure-preserving properties.
\newblock \emph{Found. Comput. Math.}, 9:\penalty0 197--219, 2009.

\bibitem[Cacuci(1981)]{Ca1981}
D.~G. Cacuci.
\newblock Sensitivity theory for nonlinear systems. {I}. {N}onlinear functional analysis approach.
\newblock \emph{J. Math. Phys.}, 22\penalty0 (12):\penalty0 2794--2802, 1981.

\bibitem[Cao et~al.(2003)Cao, Li, Petzold, and Serban]{CaLiPeSe2003}
Y.~Cao, S.~Li, L.~Petzold, and R.~Serban.
\newblock Adjoint sensitivity analysis for differential-algebraic equations: The adjoint {DAE} system and its numerical solution.
\newblock \emph{SIAM J. Sci. Comput.}, 24\penalty0 (3):\penalty0 1076--1089 (14 pages), 2003.

\bibitem[Chen et~al.(2018)Chen, Rubanova, Bettencourt, and Duvenaud]{ChRuBeDu2018}
R.~T.~Q. Chen, Y.~Rubanova, J.~Bettencourt, and D.~Duvenaud.
\newblock Neural ordinary differential equations.
\newblock In \emph{Proceedings of the 32nd International Conference on Neural Information Processing Systems}, NIPS'18, page 6572–6583, Red Hook, NY, USA, 2018. Curran Associates Inc.

\bibitem[Cohen and Welling(2016)]{CoWe2016}
T.~Cohen and M.~Welling.
\newblock Group equivariant convolutional networks.
\newblock In M.~F. Balcan and K.~Q. Weinberger, editors, \emph{Proceedings of The 33rd International Conference on Machine Learning}, volume~48 of \emph{Proceedings of Machine Learning Research}, pages 2990--2999, New York, New York, USA, 20--22 Jun 2016. PMLR.

\bibitem[Colombo et~al.(2016)Colombo, Ferraro, and Mart\'{i}n~de Diego]{Col2016}
L.~Colombo, S.~Ferraro, and D.~Mart\'{i}n~de Diego.
\newblock Geometric integrators for higher-order variational systems and their application to optimal control.
\newblock \emph{J Nonlinear Sci}, 26:\penalty0 1615–1650, 2016.
\newblock \doi{10.1007/s00332-016-9314-9}.

\bibitem[Eloe and Henderson(2016)]{ElHe2016}
Paul~W Eloe and Johnny Henderson.
\newblock \emph{Nonlinear Interpolation and Boundary Value Problems}.
\newblock Trends in Abstract and Applied Analysis: Volume 2. World Scientific, 2016.
\newblock \doi{10.1142/9877}.

\bibitem[Giles and Pierce(2000)]{GiPi2000}
M.~B. Giles and N.~A. Pierce.
\newblock An introduction to the adjoint approach to design.
\newblock \emph{Flow, Turbulence and Combustion}, 65:\penalty0 393--415, 2000.

\bibitem[Griewank(2003)]{Gr2003}
A.~Griewank.
\newblock A mathematical view of automatic differentiation.
\newblock In \emph{Acta Numer.}, volume~12, pages 321--398. Cambridge University Press, 2003.

\bibitem[Hamilton(2017)]{Ha2017}
M.~J.~D. Hamilton.
\newblock \emph{Mathematical {G}auge {T}heory}.
\newblock Springer Cham, 2017.

\bibitem[Ibragimov(2006)]{Ib2006}
N.~H. Ibragimov.
\newblock Integrating factors, adjoint equations and {L}agrangians.
\newblock \emph{Journal of Mathematical Analysis and Applications}, 318\penalty0 (2):\penalty0 742--757, 2006.

\bibitem[Ibragimov(2007)]{Ib2007}
N.~H. Ibragimov.
\newblock A new conservation theorem.
\newblock \emph{J. Math. Anal. Appl.}, 333\penalty0 (1):\penalty0 311--328, 2007.

\bibitem[Kobilarov and Marsden(2011)]{KoMa2011}
M.~B. Kobilarov and J.~E. Marsden.
\newblock Discrete geometric optimal control on {L}ie groups.
\newblock \emph{IEEE Transactions on Robotics}, 27\penalty0 (4):\penalty0 641--655, 2011.

\bibitem[Kondor and Trivedi(2018)]{KoTr2018}
R.~Kondor and S.~Trivedi.
\newblock On the generalization of equivariance and convolution in neural networks to the action of compact groups.
\newblock In Jennifer Dy and Andreas Krause, editors, \emph{Proceedings of the 35th International Conference on Machine Learning}, volume~80 of \emph{Proceedings of Machine Learning Research}, pages 2747--2755. PMLR, 10--15 Jul 2018.

\bibitem[Lee(2012)]{Le2012}
J.~M. Lee.
\newblock \emph{Introduction to Smooth Manifolds}.
\newblock Springer New York, NY, 2012.

\bibitem[Lee et~al.(2005)Lee, McClamroch, and Leok]{LeLeMc2005a}
T.~Lee, {N. H.} McClamroch, and M.~Leok.
\newblock A {L}ie group variational integrator for the attitude dynamics of a rigid body with applications to the {3D} pendulum.
\newblock \emph{Proc. IEEE Conf. on Control Applications}, pages 962--967, 2005.

\bibitem[Lee et~al.(2007{\natexlab{a}})Lee, Leok, and McClamroch]{LeLeMc2007a}
T.~Lee, M.~Leok, and {N. H.} McClamroch.
\newblock {L}ie group variational integrators for the full body problem in orbital mechanics.
\newblock \emph{Celestial Mechanics and Dynamical Astronomy}, 98\penalty0 (2):\penalty0 121--144, 2007{\natexlab{a}}.

\bibitem[Lee et~al.(2007{\natexlab{b}})Lee, Leok, and McClamroch]{LeLeMc2007b}
T.~Lee, M.~Leok, and {N. H.} McClamroch.
\newblock {L}ie group variational integrators for the full body problem.
\newblock \emph{Comput. Methods Appl. Mech. Engrg.}, 196\penalty0 (29-30):\penalty0 2907--2924, 2007{\natexlab{b}}.

\bibitem[Leok and Zhang(2011)]{LeZh2009}
M.~Leok and J.~Zhang.
\newblock Discrete {H}amiltonian variational integrators.
\newblock \emph{IMA J. Numer. Anal.}, 31\penalty0 (4):\penalty0 1497--1532, 2011.

\bibitem[Leyendecker et~al.(2024)Leyendecker, Maslovskaya, Ober-Blöbaum, de~Almagro, and Szemenyei]{Sig2024}
Sigrid Leyendecker, Sofya Maslovskaya, Sina Ober-Blöbaum, Rodrigo T. Sato~Martín de~Almagro, and Flóra~Orsolya Szemenyei.
\newblock A new {L}agrangian approach to control affine systems with a quadratic {L}agrange term.
\newblock \emph{Journal of Computational Dynamics}, 11\penalty0 (3):\penalty0 336--353, 2024.
\newblock ISSN 2158-2491.
\newblock \doi{10.3934/jcd.2024017}.

\bibitem[Li and Petzold(2003)]{LiPe2003}
S.~Li and L.~R. Petzold.
\newblock Solution adapted mesh refinement and sensitivity analysis for parabolic partial differential equation systems.
\newblock In L.~T. Biegler, M.~Heinkenschloss, O.~Ghattas, and B.~van Bloemen~Waanders, editors, \emph{Large-Scale PDE-Constrained Optimization}, pages 117--132. Springer, Berlin, 2003.

\bibitem[Ma and Rowley(2010)]{MaRo2010}
Z.~Ma and C.~W. Rowley.
\newblock {L}ie--{P}oisson integrators: A {H}amiltonian, variational approach.
\newblock \emph{Int. J. Numer. Meth. Engng.}, 82:\penalty0 1609--1644, 2010.

\bibitem[Marsden and Ratiu(1999)]{MaRa1999}
J.~E. Marsden and T.~S. Ratiu.
\newblock \emph{Introduction to Mechanics and Symmetry}.
\newblock Springer New York, NY, 1999.

\bibitem[Marsden and West(2001)]{MaWe2001}
{J. E.} Marsden and M.~West.
\newblock Discrete mechanics and variational integrators.
\newblock \emph{Acta Numer.}, 10:\penalty0 317--514, 2001.

\bibitem[Marsden et~al.(1999)Marsden, Pekarsky, and Shkoller]{MaPeSh1999}
J.~E. Marsden, S.~Pekarsky, and S.~Shkoller.
\newblock Discrete {E}uler--{P}oincar\'{e} and {L}ie--{P}oisson equations.
\newblock \emph{Nonlinearity}, 12\penalty0 (6):\penalty0 1647, 1999.

\bibitem[Matsubara et~al.(2021)Matsubara, Miyatake, and Yaguchi]{matsubara2021symplectic}
T.~Matsubara, Y.~Miyatake, and T.~Yaguchi.
\newblock Symplectic adjoint method for exact gradient of neural {ODE} with minimal memory.
\newblock In \emph{Advances in Neural Information Processing Systems}, 2021.

\bibitem[Nakahara(2003)]{Na2003}
M.~Nakahara.
\newblock \emph{Geometry, {T}opology, and {P}hysics}.
\newblock CRC Press, second edition, 2003.

\bibitem[Nguyen et~al.(2016)Nguyen, Georges, and Besan\c{c}on]{NgGeBe2016}
V.~T. Nguyen, D.~Georges, and G.~Besan\c{c}on.
\newblock State and parameter estimation in 1-{D} hyperbolic {PDEs} based on an adjoint method.
\newblock \emph{Automatica}, 67\penalty0 (C):\penalty0 185--191, May 2016.
\newblock ISSN 0005-1098.

\bibitem[Pierce and Giles(2000)]{PiGi2000}
N.~A. Pierce and M.~B. Giles.
\newblock Adjoint recovery of superconvergent functionals from {PDE} approximations.
\newblock \emph{SIAM Rev.}, 42\penalty0 (2):\penalty0 247--264, 2000.

\bibitem[Ross(2005)]{Ro2005}
I.~M. Ross.
\newblock A roadmap for optimal control: The right way to commute.
\newblock \emph{Ann. NY Acad. Sci.}, 1065\penalty0 (1):\penalty0 210--231, 2005.

\bibitem[Ross and Fahroo(2001)]{RoFa2001}
M.~Ross and F.~Fahroo.
\newblock A pseudospectral transformation of the convectors of optimal control systems.
\newblock \emph{IFAC Proc. Ser.}, 34\penalty0 (13):\penalty0 543--548, 2001.

\bibitem[Sanz-Serna(2016)]{Sa2016}
J.~M. Sanz-Serna.
\newblock Symplectic {R}unge--{K}utta schemes for adjoint equations, automatic differentiation, optimal control, and more.
\newblock \emph{SIAM Review}, 58\penalty0 (1):\penalty0 3--33, 2016.

\bibitem[Sirkes and Tziperman(1997)]{SoTz1997}
Z.~Sirkes and E.~Tziperman.
\newblock Finite difference of adjoint or adjoint of finite difference?
\newblock \emph{Mon. Weather Rev.}, 125\penalty0 (12):\penalty0 3373--3378, 1997.

\bibitem[Tran and Leok(2022)]{TrLe2022}
B.~Tran and M.~Leok.
\newblock Multisymplectic {H}amiltonian variational integrators.
\newblock \emph{International Journal of Computer Mathematics (Special Issue on Geometric Numerical Integration, Twenty-Five Years Later)}, 99\penalty0 (1):\penalty0 113--157, 2022.

\bibitem[Tran and Leok(2024)]{TrLe2022adj}
B.~K. Tran and M.~Leok.
\newblock Geometric methods for adjoint systems.
\newblock \emph{J Nonlinear Sci}, 34\penalty0 (25), 2024.
\newblock \doi{10.1007/s00332-023-09999-7}.

\bibitem[Tran et~al.(2024)Tran, Southworth, and Leok]{TrSoLe2024}
B.~K. Tran, B.~S. Southworth, and M.~Leok.
\newblock On properties of adjoint systems for evolutionary pdes.
\newblock \emph{J Nonlinear Sci}, 34\penalty0 (95), 2024.
\newblock \doi{10.1007/s00332-024-10071-1}.

\bibitem[Wang et~al.(2012)Wang, Duraisamy, Alonso, and Iaccarino]{WaDuAlIa2012}
Q.~Wang, K.~Duraisamy, J.~J. Alonso, and G.~Iaccarino.
\newblock Risk assessment of scramjet unstart using adjoint-based sampling.
\newblock \emph{AIAA J.}, 50\penalty0 (3):\penalty0 581--592, 2012.

\end{thebibliography}

\end{document}